\input amstex
\documentstyle{amsppt}
\magnification = \magstep1
\topmatter
\title Curves in the complement of a smooth plane cubic 
whose normalizations are $\Bbb A^1$
\endtitle
\author Nobuyoshi Takahashi 
\footnote{Research Fellow of the Japan Society for the Promotion of Sciences. 
This work was supported by the Grant-in-Aid for Scientific Research  of 
Monbusho(Ministry of Education, Science and Culture, Japan).}
\endauthor
\rightheadtext{$\Bbb A^1$'s in $\Bbb P^3\backslash(\text{cubic})$}
\leftheadtext{Nobuyoshi Takahashi}
\address Department of Mathematical Sciences, University of Tokyo, Komaba, 
Meguro, Tokyo, 153, Japan \endaddress
\email nn57016\@hongo.ecc.u-tokyo.ac.jp \endemail
\curraddr Ecole Normale Sup\'erieure, 48, boulevard Jourdan, 
75690 Paris Cedex 14, France \endcurraddr
\email takahasi\@clipper.ens.fr \endemail
\abstract
 For a smooth plane cubic $B$,
we count curves $C$ of degree $d$
such that the normalizations of $C\backslash B$
are isomorphic to $\Bbb A^1$, for $d\leq7$
(for $d=7$ under some assumption).
We also count plane rational quartic curves
intersecting $B$ at only one point. 
\endabstract
\endtopmatter
\document
\head \S0 Introduction \endhead

 Recently, the space of rational curves in a manifold 
is attracting attention in connection with physics. 
The number of rational curves of given degree in a Calabi-Yau threefold 
is studied, inspired by mirror symmetry(see [Ms], [Kn]), 
and the degree of the variety of rational curves in a fixed homology class 
in Fano manifolds 
is computed using the associativity of the quantum product([KM], [RT]). 

 In this paper, for a smooth rational surface $X$ (over $\Bbb C$)
and its smooth anti-canonical divisor $B$,
we count curves $C$ in a fixed divisor class in $X$
such that the normalizations of $C\backslash B$ are isomorphic to $\Bbb A^1$.
In this section, we will refer to them, or to $C\backslash B$,
as affine lines in $X\backslash B$ for simplicity. 

 On the one hand,
this can be considered as an analog of counting rational curves
on a K3 surface or a Calabi-Yau threefold,
and it is hoped that they enjoy some good property 
similar to ``mirror phenomena'' etc.
(while a rational curve is considered as the locus swept by a closed string, 
an affine line might be seen as that swept by an infinitely long string),
although I don't know what it should be. 
On the other hand, this problem can be seen as a special case 
of counting rational curves in rational surfaces. 
In fact, our computations in this paper use 
the degree of the variety of rational curves in rational surfaces([CM]). 

 The problem has also an arithmetic application([B], [V]).
If $K$ is a number field,
$S$ a finite set of places of $K$ containing all infinite ones,
$\Cal O_S$ the ring of $S$-integers
and $B_S\subset\Bbb P^2_{\Cal O_S}$ a cubic curve smooth over
$\text{Spec}\Cal O_S$,
then an $S$-integral point of $\Bbb P^2_{\Cal O_S}\backslash B_S$ is,
by definition, its $\Cal O_S$-valued point.
If $C_S$ is a curve in $\Bbb P^2_{\Cal O_S}\backslash B_S$,
it follows from the Siegel-Mahler theorem
that there exists a finite extension $K'$ of $K$
and a finite set $S'$ of places of $K'$ containing all infinite ones
such that there are infinitely many $S'$-integral points contained in
$C_S\otimes \Cal O_{S'}$, 
if and only if $C_S\otimes\Bbb C$ is isomorphic to $\Bbb A^1$
or $\Bbb G_m$. 
The number of $S'$-integral points of $C_S\otimes \Cal O_{S'}$
whose heights are less than $N$
is of order $N^a(\log N)^b$ with $a>0, b\geq0$ 
in the former case, 
while it is of order $(\log N)^b$ with $b>0$
in the latter case.
Thus affine lines in $\Bbb P^2\backslash(B_S\otimes\Bbb C)$ 
give many solutions to the corresponding Diophantine equation. 

 In particular, let $X=\Bbb P^2$.
If $C$ is of degree $d$,
$C\cap B$ is seen to consist of a point of order dividing $3d$
for the addition on $B$ with a flex as zero.
We consider mainly the number $n_d$, 
counted with appropriately defined multiplicity, 
of such curves with $C\cap B$ a fixed point of order $3d$ w.r.t. some flex. 
We obtain the answer for $d\leq6$, 
and for $d=7$ under some assumption. 

 A useful fact is that $\Bbb P^2$ has a triple cover $Y$
totally ramified over $B$ and unramified elsewhere:
if $B$ is defined by $F(X,Y,Z)=0$,
$Y$ can be written as $F=W^3$ in $\Bbb P^3$.
In \S2, we obtain $n_4=16$ using this(Theorem 2.1). 
We also prove that $n_d>0$ holds for any $d$ 
and that, if $B_S$ is a smooth cubic in $\Bbb P^2_{\Cal O_S}$, 
where $S$ is a finite set of places of a number field $K$
containing all infinite ones,
there exist a finite extension $K'$ of $K$
and a finite set $S'$ of places of $K'$, containing all infinite ones,
such that, for any $d$, there exists a morphism
$\Bbb A^1_{\Cal O_{S'}}\rightarrow\Bbb P^2_{\Cal O_{S'}}
\backslash(B_S\otimes\Cal O_{S'})$ 
over $\Cal O_{S'}$ 
which is birational onto the image (on the generic fiber) 
and is of degree larger than $d$ (Proposition 2.3). 

 For $d=5$ and $d=6$, using the same trick,
we reduce the problem to counting affine lines in different divisor classes
in an affine cubic surface.
In \S3, as a preparation for the computation of these numbers,
we construct a space 
representing affine lines and some of their degeneration in $X\backslash B$,
where $X$ is a smooth rational surface 
and $B$ is a smooth anti-canonical divisor, 
and find relations between the numbers of affine lines in $X\backslash B$ 
for different $X$ and $B$(Theorem 3.8). 
Then, we obtain $n_5=113$ and $n_6=948$(Corollary 4.4)
from the number of rational plane curves of degree $\leq4$ 
which have given multiplicities at given general points. 

 For $d=7$, it is necessary to know the number $x$ of rational sextic curves
which have multiplicity 2 at 8 general points
and pass through another general point,
but it seems that this number cannot be obtained
in the same way as the number of, say, rational sextic curves
which have multiplicity 2 at 7 general points
and pass through other 3 general points.
By calculating $n_6$ without using the triple cover,
we see $x=90$ under some assumption. 
Then we obtain $n_7=8974$(Corollary 4.5). 

 Using Theorem 3.8, 
we also compute the number of quartic curves 
which meet $B$ at only one point of order 12(Corollary 4.3).

 Calculations for the case $d=3$ are also in [R], [V] and [X].
[V] and [X] are rather interested in
the number of rational curves having only one point on the cubic curve,
and [V] also contains some explicit equations of such curves for $d=4$.
[X] proves a theorem concerning finiteness of curves in surfaces 
which has only one point of intersection with a fixed curve, 
which is more general than Proposition 1.1(with $n=2$).

{\bf Acknowledgements.} 
 The author would like to express his sincere gratitude 
to Professor Y.~Kawamata for helpful advice and encouragement. 
He is very grateful to Professor Y.~Tschinkel 
for instructing him in the arithmetic background of this problem. 
He would like to thank Ecole Normale Sup\'erieure, 
where this work was done, for hospitality. 

\head \S1 Case $d=1,2 \text{ and } 3$ \endhead

 In this paper, we work over $\Bbb C$ unless otherwise stated. 

 Let $B$ be a smooth cubic curve in $X=\Bbb P^2$.
We want to count irreducible reduced curves $C$ of degree $d$ in $X$
such that the normalizations of $C\backslash B$ are isomorphic to $\Bbb A^1$
(or, equivalently, rational curves $C$ having only one place on $B$), 
for, according to a rough dimension counting argument, 
the expected dimension of the variety of such curves is 0. 
In fact, we have: 

\proclaim{Proposition 1.1}
 Let $X$ be a smooth projective variety 
and $B\subset X$ a smooth hypersurface
such that $\vert K_X+B\vert\not=\phi$.
Then, for any homology class $\beta$ on $X$,
the union of rational curves in the class $\beta$
having only one point of intersection with $B$ is 
contained in a proper closed subset of $X$ 
(in the Zariski topology). 
\endproclaim
\demo{Proof}
 Assuming the contrary, we have a smooth variety $S$ of dimension $n-1$,
disjoint sections $S_i$ of $p_1:Y = S \times \Bbb P^1 \rightarrow S$
and a dominant morphism $f:Y\rightarrow X$
such that $f^{*}(B)=\sum a_i S_i$, where $a_i$ are positive integers,
and that 
$f\circ (p_1\vert_{S_i})^{-1}=f\circ (p_1\vert_{S_j})^{-1}$ for all $i,j$.
Denote this last morphism by $g:S\rightarrow B$. 

 On $X$, there exists a nonzero $n$-form $\omega$, where $n=\dim X$,
which is regular outside $B$ and has at most logarithmic pole along $B$.
Then the pullback $\omega_Y = f^* \omega$(as a differential form) 
is a nonzero $n$-form on $Y$ 
which has only logarithmic pole along $\cup S_i$, 
and we have $\text{res}_{S_i}\omega_Y = a_i f^{*}\text{res}_B \omega$.
If we take a non-vanishing $(n-1)$-form $\omega_S$ on $S$,
shrinking $S$ if necessary,
we can write $\omega_Y$ as $\omega_{Y/S}.p_1^{*}\omega_S$
for some relative differential 1-form $\omega_{Y/S}$.
Now restricting $\omega_{Y/S}$ to a general fiber of $p_1$,
we obtain a nonzero 1-form on $\Bbb P^1$ with at most logarithmic poles 
whose residues are positive integer times some complex number, 
which is impossible. \qed 
\enddemo 

 The computations for $d= 1 \text{ and } 2$ are easy. First, we have 

\proclaim{Lemma 1.2}
Let $B$ be a smooth plane cubic curve
and fix a flex $O$ of $B$.
Then, for a curve $C$ of degree $d$ meeting $B$ at only one point, 
$C\vert _B = 3dP$, where $P$ is a point of order dividing $3d$ 
with respect to the addition on $B$ with origin $O$. \qed
\endproclaim

 So we make the following: 

\definition{Definition 1.3}
 Let $X$ be a rational surface, $B$ a smooth anti-canonical divisor on $X$
and $P$ a point on $B$.
We say a curve $C$ on $X$ satisfies (*) for $(B, P)$
if $C$ is irreducible and reduced,
the normalization of $C\backslash B$ is isomorphic to $\Bbb A^1$
or $\Bbb P^1$ 
and $C\cap B\subseteq\{P\}$. 
\enddefinition

 Using the notation in the lemma above, we have 

\proclaim{Proposition 1.4}
 For a flex $P$, there is precisely one line 
and no conic satisfying (*) for $(B, P)$. 
For a point $P$ of order 6, 
there is exactly one conic satisfying (*) for $(B, P)$. 
\endproclaim

\demo{Proof}
 For a point $P$ of order dividing $3d$, 
$3dP \sim dH\vert_B$ and therefore the long exact sequence 
associated to the exact sequence 

$$
  0 \rightarrow \Cal O_X(d-3) \rightarrow \Cal O_X(d) \rightarrow 
                \Cal O_B(dH\vert_B) \rightarrow 0 
$$
shows, for $d=1 \text{ and } 2$, that there is unique curve of degree $d$ 
intersecting with $B$ only at $P$. 
For $d=1$, this curve is a line and satisfies (*). 
For $d=2$, if $P$ is a flex, then this curve is 2 times line 
and so does not satisfy (*). 
If $P$ is not a flex, then this curve 
cannot be reducible or non-reduced, as it would mean that 
there is a line having triple intersection with $B$ at $P$. \qed
\enddemo

 The case $d=3$ is a little more difficult. 

\proclaim{Proposition 1.5}
 For a point $P$ of order 9, 
there are 3 cubic curves satisfying (*) for $(B, P)$, 
and they have only nodes. 
For a flex $P$, if $B$ is not isomorphic to the curve $B_0$ 
defined by $y^2=x^3-1$, 
there are two cubic curves satisfying (*) for $(B, P)$, 
and they have only nodes. 
If $B$ is isomorphic to $B_0$, then there is unique cubic
satisfying (*) for $(B, P)$, and it has a cusp outside $B$. 
\endproclaim
\demo{Proof}
 By the same exact sequence as in the proof of the proposition above,
we see that the cubic curves
having intersection of multiplicity 9 at $P$ with $B$
together with $B$ itself form a pencil $\Lambda$. 
Blowing up 9 times the points on the proper transform of $B$ over $P$, 
we obtain a resolution of $\Lambda$, giving an elliptic fibration 
$f:X' \rightarrow \Bbb P^1$. 

 We see that there is just one member $D$ of $\Lambda$ 
which is singular at $P$, 
for being singular at $P$ is a linear condition, 
and that the proper transform of $D$ and the exceptional curves 
other than the last one form a set theoretic fiber $F$ of $f$. 

 If $P$ is not a flex, there is no reducible or non-reduced member 
in $\Lambda$ by the same argument as in Proposition 1.4. 
Furthermore, $D$ cannot have a cusp at $P$, for a cusp cannot have 
a 9-ple intersection with a smooth curve. 
Therefore $D$ has a node at $P$, and $F$ is a cycle 
consisting of 9 smooth rational curves. 
Assume that $\Lambda$ contains a cuspidal cubic $E$. 
Then for the addition on $E\backslash\{\text{cusp}\}$ 
with the flex as origin, $P$ is a 9-torsion, 
and therefore $P$ is the flex of $E$, 
which is impossible as $P$ is not a flex of $B$ 
and $B$ and $E$ have 9-ple intersection at $P$. 
Now the Euler number of $X'$ is 12, and that of $F$ is 9. 
Consequently, the number of fibers of $f$ isomorphic to a nodal cubic is 3, 
which is the number of curves satisfying (*) for $(B, P)$. 

 If $P$ is a flex, $D$ is 3 times the tangent line 
and the Euler number of $F$ is 10. 
If $\Lambda$ contains a cuspidal cubic $E$, 
then $P$ is a flex of $E$. 
But then, writing as $E:Y^2 Z = X^3$ and $P(0:1:0)$, 
we see that $B$ must be of the form $Y^2 Z = X^3 + aZ^3$ for $a\in \Bbb C$, 
i.e. $B$ is isomorphic to $B_0$. 
In this case, since the Euler number of a cuspidal cubic is 2, 
$\Lambda$ contains one cuspidal cubic and no nodal cubic. 
If $B$ is not isomorphic to $B_0$, 
$\Lambda$ contains 2 nodal cubics. \qed
\enddemo

\head \S2 Case $d=4$ \endhead

 In this section, we compute the number for $d=4$ and $P$ ``primitive''. 

\proclaim{Theorem 2.1}
 Let $B$ be a general cubic and $P$ a point of order 12. 
Then any curve of degree 4 
satisfying (*) for $(B, P)$ has only nodes 
and there are 16 such curves. 
\endproclaim

\demo{Proof}
 Let $C$ be a curve satisfying (*) for $(B, P)$.
Then, since $C$ and $B$ intersects at $P$ with multiplicity 12,
$C$ is smooth at $P$.
Assume that $C$ has a cusp $Q$ and two nodes $R$, $S$, for example.
Denote the normalization by $f:\Bbb P^1\rightarrow C$ and
the curve obtained by resolving the singularities $R$ and $S$
(resp. $S$ and $Q$, $Q$ and $R$) of $C$ by $C_1$(resp. $C_2$, $C_3$).
Take a coordinate $u$ on $\Bbb P^1$ such that
$f^{-1}(Q)=\{u=\infty\}$
and that $f^{-1}(R)=\{u=0,1\}$.
Let $f^{-1}(S)=\{u=a,b\}$ and $f^{-1}(P)=\{u=c\}$.
We have $12P\sim3H\vert_C$.
Considering the lines through $R$ and $S$,
$S$ and $Q$, $Q$ and $R$, we have
$3([0]+[1]+[a]+[b])\sim 12[c]$ on $C_1$,
$3([a]+[b]+2[\infty])\sim 12[c]$ on $C_2$
and $3([0]+[1]+2[\infty])\sim 12[c]$ on $C_3$. 
Thus $3(a+b+1) = 12c$,
$(ab/(a-1)(b-1))^3= (c/(c-1))^{12}$, and
$(a(a-1)/b(b-1))^3= ((c-a)/(c-b))^{12}$.
Eliminating $c$ and computing the greatest common divisor of 
the resulting polynomials, 
we see that any positive dimensional component
of the solution of the equation above 
satisfies $[0]+[1]+[a]+[b]\sim 4[c]$ on $C_1$,
$[a]+[b]+2[\infty]\sim 4[c]$ on $C_2$
and $[0]+[1]+2[\infty]\sim 4[c]$ on $C_3$.
Then $4[c]$ is linearly equivalent to the pullback of $H$ on $C_i$,
and therefore on $C$.
Writing down a long exact sequence, we see that there exists a line 
intersecting $C$ with multiplicity 4 at $P$.
This line intersects with $B$ at $P$ with multiplicity 4,
which is absurd. 
Thus there are only finitely many possibilities
for $C$ and $P$ modulo projective equivalence.
Since $B$ and $C$ intersect with multiplicity 12 at $P$,
$B$ is determined uniquely by $C$ and $P$.
Therefore, for a general $B$, there does not exist such a curve $C$.
The proof that $C$ cannot have even worse singularity 
is similar and easier. 

 Let $\pi:Y\to X$ be the triple cover totally branched along $B$,
defined by $t^3=s$ in the total space
of the line bundle associated to $\Cal O(1)$,
where $t$ is a local coordinate of that bundle
and $s$ is a section of $\Cal O(3)$ with $B=(s=0)$.
The line bundle can be embedded in $\Bbb P^3$,
and then $Y$ is represented as a cubic surface:
in fact, if $B$ is defined by $F(X, Y, Z)=0$,
$Y$ can be written as $F(X, Y, Z)=W^3$. 
Let $B_Y$ be the reduced inverse image of $B$,
and then it is easy to see that $B_Y$ is a plane section
and maps isomorphically onto $B$.
We identify points on $B$ and $B_Y$. 

 Let $C$ be a curve satisfying (*) for $(B, P)$.
Since the covering of $C$ by $\pi^*(C)$ is unbranched outside $P$,
$\pi^*(C)$ is composed of 3 curves satisfying (*) for $(B_Y, P)$,
each mapping birationally to $C$ so that any local analytic branch
maps isomorphically onto the corresponding branch.
Furthermore, since they are permuted cyclically
by the action of $\text{Gal}(Y/X)$,
they have the same degree, i.e. 4.
On the other hand, if $C$ is a rational curve of degree 4 on $Y$
which satisfies (*) for $(B_Y, P)$
then $\pi(C)$ is a curve satisfying (*) for $(B, P)$,
and therefore $\pi:C\rightarrow \pi(C)$ is birational
and maps any analytic branch isomorphically, from what we have seen above
(and therefore $C$ has only nodes). 
Thus it suffices to prove that 
there are 48 curves of degree 4 on $Y$ which satisfy (*) for $(B_Y, P)$. 

 To prove this, let $g:Y\to Z\cong \Bbb P^2$ be a blow-down.
Since the 27 lines on $Y$ are the components of the inverse images
of the flex tangent lines to $B$,
$g$ is the blow-down of 6 disjoint lines $E_1, \dots, E_6$
intersecting with $B_Y$ at different 3-torsions $P_1, \dots, P_6$
(with a flex $O$ on $B\subset X$ as zero).
Let $Q_1, Q_2, Q_3$ be the other 3-torsions,
$B_Z = g_*B_Y$, which is a cubic curve,
and $O_Z$ a flex of $B_Z$.
None of 3 lines on $Y$ through $Q_1$ is contracted by $g$, and
the sum of the degrees of their images is easily seen to be 3.
Therefore, they are lines, and any of them meets $B_Z$
at $Q_1, P_a \text{ and } P_b$ for some $a\not=b$. 
Thus we have $3O_Z \sim Q_1 + P_a + P_b$
and $O_Z$ is a 9-torsion relative to $O$.
Now $B$ is a group isomorphic to $(\Bbb R/\Bbb Z)^2$.
By changing $O$ if necessary, we can find such an identification with
$Q_1=(0, 2/3), Q_2=(1/3, 2/3)$ and $Q_3=(2/3, 2/3)$ or $Q_3=(2/3, 1/3)$.
As any pairwise distinct three of $P_i$ cannot be on a line in $Z$,
for $B_Y$ is ample, 
we see that $Q_3=(2/3, 2/3)$ and $3O_Z=(1/3, 0)$ or $(2/3, 0)$.
We may assume $O_Z=(1/9, 0)$. 

 Let $C$ be a curve as above.
If $C\sim eH-\sum a_iE_i$, where $H$ denotes the pullback of a line in $Z$, 
we have $0\leq a_i=E_i.C\leq \pi^*(\text{line}).C=4$
and $e=\deg g_*C = (B_Y+E_1+\dots+E_6).C/3=(4+\sum a_i)/3$,
and therefore $2\leq d\leq9$. 

 Using furthermore $p_a(C)=(e-1)(e-2)/2-\sum a_i(a_i-1)/2\geq0$,
we have the following list for $e$, $p_a$ and unordered sequences 
$[a]=[a_1, \dots, a_6]$:

$$\vbox{\halign{&\hfil#\hfil \cr
  $e$ & $p_a$ & $[a]$ & \quad & $e$ & $p_a$ & $[a]$& \quad & 
     $e$ & $p_a$ & $[a]$ \cr
   \noalign{\hrule} \noalign{\smallskip}
  2 & 0 & $[0,0,0,0,1,1]$ & \quad & 4 & 0 & $[0,1,1,2,2,2]$ & \quad & 
     5 & 0 & $[1,1,2,2,2,3]$   \cr
  3 & 0 & $[0,0,1,1,1,2]$ & \quad &   & 0 & $[1,1,1,1,1,3]$ & \quad & 
       & 1 & $[1,2,2,2,2,2]$   \cr
    & 1 & $[0,1,1,1,1,1]$ & \quad &   & 1 & $[1,1,1,1,2,2]$ & \quad &  
     6 & 0 & \,$[2,2,2,2,3,3]$.\cr
}}$$

 It can be seen that given such a curve $C$
we can make another choice of $Z$, say $Z'$,
such that $e'=2, [a']=[0,0,0,0,1,1]$ for case $p_a=0$
or that $e'=3, [a']=[0,1,1,1,1,1]$ for case $p_a=1$,
where we denote the degrees etc. for $Z'$ by symbols with $'$:
for example, in the case $e=6$,
it suffices to apply two successive quadratic transformations, 
first with centers $P_k, P_l, P_m$ with $a_k=a_l=3$ 
and then with the other three points as the centers. 

 For the case $p_a=0$ above,
there is one (possibly reducible or nonreduced) conic on $Z'$
having 4-fold intersection at $P$
and simple intersections at $P_k$ and $P_l$ with $B_{Z'}$,
where $a'_k=a'_l=1, k\not=l$,
if and only if $4P+P_k+P_l \sim 6O_{Z'}$.
By $eH-\sum a_iE_i\sim2H'-E'_k-E'_l$, 
this is equivalent to say $4(P-O_Z) + \sum a_i(P_i-O_Z) \sim 0$. 
Now we have the following: \enddemo

\proclaim{Lemma 2.2}
 Let $D\in\vert eH-\sum a_iE_i\vert$ be an effective divisor on $Y$
such that $D\vert_{B_Y}=dQ$,
where $d:=3e-\sum a_i$ and $Q\not\in\cup E_i$ is a point of order $3d$ 
w.r.t. a flex $O\in B\subset X$. 
Then $C$ is irreducible and reduced. 
\endproclaim
\demo{Proof}
 If $C$ is reducible or reduced,
we have $d'P+\sum a'_iP_i\sim e'H$ with $0<d'<d$, since $B_Y$ is ample.
This implies $3d'(P-O)\sim0$, a contradiction. \qed
\enddemo

Therefore this conic gives
a smooth rational quartic curve on $Y$
which satisfies (*) for $(B_Y, P)$ 
and is linearly equivalent to $eH-\sum a_iE_i$. 

 On the other hand, for the case $p_a=1$,
there is a linear pencil on $Z'$
consisting of $B_{Z'}$ and cubics on $Z'$ 
having 4-fold intersection at $P$
and simple intersection at $P_1, \dots, \breve P_k, \dots, P_6$
with $B_{Z'}$, where $a'_k=0$,
if and only if $4P+P_1+\dots+\breve P_k+\dots+P_6 \sim 9O_{Z'}$,
which is again equivalent to $4(P-O_Z) + \sum a_i(P_i-O_Z) \sim 0$.
By Lemma 2.2,
all members of the pencil, including $B_{Z'}$ itself, 
are irreducible and reduced. 
As in Proposition 1.5, 
we see that the member which is singular at $P$ is nodal 
and therefore that there are 8 nodal cubics in the pencil 
with only one place on $P$, 
giving nodal rational quartic curves on $Y$
which satisfy (*) for $(B_Y, P)$ 
and are linearly equivalent to $eH-\sum a_iE_i$. 

 By computation(with a computer), we see that the number of ordered sequences
$(a_1, a_2, a_3, a_4, a_5, a_6)$ 
such that $4(P-O_Z) + \sum a_i(P_i-O_Z) \sim 0$ is: 

\medpagebreak
\noindent for $p_a=0$, 

$$\vbox{\halign{&\hfil#\hfil \cr
    $e$  & $[a]$ & if $4P=(x, 0)$ or $(x, 1/3)$ \quad & if $4P=(x, 2/3)$ \cr
     \noalign{\hrule} \noalign{\smallskip}
     2   & $[0,0,0,0,1,1]$ & 1 & 3 \cr
     3   & $[0,0,1,1,1,2]$ & 7 & 6 \cr
     4   & $[0,1,1,2,2,2]$ & 7 & 6 \cr
         & $[1,1,1,1,1,3]$ & 1 & 0 \cr
     5   & $[1,1,2,2,2,3]$ & 7 & 6 \cr
     6   & $[2,2,2,2,3,3]$ & 1 & \;3, \cr
}}$$
for $p_a=1$, 

$$\vbox{\halign{&\hfil#\hfil \cr
    $e$  & $[a]$ & if $4P=(x, 0)$ or $(x, 1/3)$ \quad & if $4P=(x, 2/3)$ \cr
     \noalign{\hrule} \noalign{\smallskip}
     3   & $[0,1,1,1,1,1]$ & 1 & 0 \cr
     4   & $[1,1,1,1,2,2]$ & 1 & 3 \cr
     5   & $[1,2,2,2,2,2]$ & 1 & \;0. \cr
}}$$

 Consequently, for any point $P$ of order 12, 
there are 48 quartic rational curves on $Y$
satisfying (*) for $(B_Y, P)$, as desired. \qed

\remark{Remark}
 It may be more natural to study, from the start, the curves
in $Y\backslash B_Y$,
where $Y$ is a general cubic surface 
and $B_Y$ is a general plane section, 
whose normalizations are $\Bbb A^1$. 
\endremark

\proclaim{Proposition 2.3}
 (1) Let $B$ be a smooth cubic curve and $d$ a positive integer.
Then there exists a rational curve $C$ in $\Bbb P^2$ of degree $d$
such that the normalization of $C\backslash B$ is
isomorphic to $\Bbb A^1$
and that the normalization map $\Bbb P^1\rightarrow C$ is immersive. 

 (2) Let $K$ be a number field,
$S$ a finite set of places of $K$ containing all infinite ones
and $B_S$ a divisor of degree 3 in $\Bbb P^2_{\Cal O_S}$
which is generically smooth.
Then there exist a finite extension $K'$ of $K$
and a finite set $S'$ of places of $K'$, containing all infinite ones,
which satisfy the following conditions: 

 (a) $\Cal O_{S'}\supset\Cal O_S$. 

 (b) Let $B_{S'}:=B_S\otimes\Cal O_{S'}$. 
Then, for any $d$, there exists a morphism 
$\Bbb A^1_{\Cal O_{S'}}\rightarrow\Bbb P^2_{\Cal O_{S'}}\backslash B_{S'}$
over $\Cal O_{S'}$
which is immersive, birational onto the image
and of degree larger than $d$ on the generic fiber. 
\endproclaim

\demo{Proof}
 We use the notations in the proof of Theorem 2.1. 

 (1) If $d=2k-1$, let $D:=kH-(k-1)E_1-E_2-E_4$ and $P:=((k-1)/3d, 1/3d)$,
and if $d=2k$, let $D:=kH-(k-1)E_1-E_4$ and $P:=(k/3d,-1/3d)$.
We may assume that $k\geq2$.
Then there exists an effective divisor $C\in\vert D\vert$ on $Y$
with $C\vert{B_Y}=dP$:
in fact, curves $C$ in $Z$ of degree $k$ such that
$C\vert_{B_Z}=dP+(k-1)P_1+P_2+P_4$
(resp. $C\vert_{B_Z}=dP+(k-1)P_1+P_4$) or $C\supset B$
form a linear system of dimension $k-2$,
and there exists such a curve with multiplicity $k-1$ at $P_1$.
Clearly this cannot contain $B_Z$. 

 Then $C$ is irreducible and reduced by Lemma 2.2, 
and therefore it is a smooth rational curve. 

 (2) Take $K'$ and $S'$ such that: 

 ($\alpha$) (a) holds, and $B_{S'}$ is smooth over $\text{Spec}\Cal O_{S'}$, 

 ($\beta$) one of the flexes of $B_{S'}$ is defined over $\Cal O_{S'}$, 

 ($\gamma$) and there exists an open immersion 
$h:\Bbb P^2_{\Cal O_{S'}}\backslash B_{S'}\rightarrow Y_{S'}$ 
whose inverse is the contraction of six lines, 
where $\pi:Y_{S'}\rightarrow\Bbb P^2_{\Cal O_{S'}}$ 
is the triple cover determined by $B_S'$. 

 Taking a flex tangent line,
we have a morphism of degree 1 as in the assertion.
If we have such a morphism $f$ of degree $d$,
then $\pi\circ h\circ f$ is also a morphism as in the assertion. 
It is of degree $3d$ if the image of the point at infinity by $f$
is not one of the fundamental points of $\pi^{-1}$,
$3d-1$ if it is. \qed
\enddemo

\head \S3 A space representing affine lines \endhead

 Let $X$ be a smooth rational surface and $B$ a smooth anti-canonical divisor.
In this section,
we obtain relations between the number of curves in $X\backslash B$
whose normalizations are $\Bbb A^1$ and similar numbers
on the blow-up of $X$ at a point on $B$. 

\medpagebreak

 Let $d$ be a positive integer. 
Take two copies of $\Bbb P^1$ and denote the second by $S^{(d)}$. 
In this section, $t\in\Bbb C\cup\{\infty\}$ 
always denotes a inhomogeneous coordinate of $S^{(d)}$ 
and $s:=1/t$ the coordinate near $\infty$. 
Consider $C':=\Bbb P^1\times S^{(d)}$ with the second projection $p'$,
four sections $\sigma'_1:=\{0\}\times S^{(d)}, 
\sigma'_2:=\{1\}\times S^{(d)}, \sigma'_4:=\{\infty\}\times S^{(d)}$ and
$\sigma'_3:=(\text{diagonal})$.
Blow up $C'$ successively $d$ times
at the intersection of the proper transform of $\sigma'_3$
and the fiber over $\infty$, and then blow down the exceptional curves
other than the last one.
Call this surface $C^{(d)}$, 
with the projection $p^{(d)}:C^{(d)}\rightarrow S^{(d)}$
and sections $\sigma_1, \sigma_2, \sigma_3$ and $\sigma_4$.
We write the fibers as $C^{(d)}_t$, $\sigma_{1,t}$, etc.
A point in this construction is that $D^{(d)}:=d\sigma_4$ is Cartier. 

 On the other hand,
let $X$ be a projective variety, $B$ a Cartier divisor on $X$,
$\beta$ a homology class of degree 2
and $H_1$, $H_2$ and $H_3$ Cartier divisors.
Then let $\bar{M}$ be the closed subscheme
of $\text{Mor}_{S^{(d)}}(C^{(d)}, X\times S^{(d)})$
parameterizing morphisms $C^{(d)}_t\rightarrow X$ 
with images in the class $\beta$ 
by which the pullback of $H_i$ contains $\sigma_{i,t}$ for $i=1, 2, 3$ 
and that of $B$ contains $D^{(d)}_t$. 
Denote by $M(H_1, H_2, H_3)$ the open subscheme of $\bar{M}$
whose geometric points represent birational morphisms
from $C^{(d)}_t$ to its images,
where $\sigma_{i,t}$ ($i=1,2,3,4$) are distinct points(i.e. $t\not=0, 1$),
by which the pullback of $H_i$ is exactly $\sigma_{i,t}$ 
in a neighborhood of $\sigma_{i,t}$ for $i=1, 2, 3$
and that of $B$ is $D^{(d)}_t$ in a neighborhood of
$\text{Supp}D^{(d)}_t$.

\proclaim{Proposition 3.1}
 (1) Let $\tilde{M}$ be a disjoint union of $M(H_1, H_2, H_3)$ for
a finite set of $(H_1, H_2, H_3)$, 
$\tilde{p}:\tilde{C}\rightarrow \tilde{M}$, 
$\tilde{f}:\tilde{C}\rightarrow X$ 
and $\tilde{D}$ the representing family. 
Define $I$ to be 
$\text{Isom}_{\tilde{M}\times\tilde{M}} 
((\tilde{C}, \tilde{D}, \tilde{f})\times 
\tilde{M}, \tilde{M}\times(\tilde{C}, \tilde{D}, \tilde{f}))$,
i.e. the scheme parameterizing
isomorphisms of different fibers of $\tilde{p}$
commuting with the restrictions of $\tilde{f}$
and inducing isomorphisms of the restrictions of $\tilde{D}$.
Then the structure morphism $\pi:I\rightarrow \tilde{M}\times\tilde{M}$
is a closed immersion
and defines an \'etale equivalence relation on $\tilde{M}$.
Therefore it defines a quotient $M$ of $\tilde{M}$
as a separated algebraic space,
as well as quotients $C$, $D$ and $f$ of $\tilde{C}, \tilde{D}, \tilde{f}$.
For some finite collection of $H_i$'s,
$M$ contains any $M(H_1, H_2, H_3)$
as an \'etale open set.

 (2) Similarly, for $d=(B.\beta)=0$, there exists an algebraic space $M$
parameterizing birational morphisms $\Bbb P^1\rightarrow X\backslash B$ 
whose images are in the class $\beta$ 
modulo automorphisms of $\Bbb P^1$. 
\endproclaim

\demo{Proof}
 (1) First, we show that $\pi$ is a closed immersion.
It is easy to see that any geometric fiber of $\pi$
is empty or consists of one reduced point, since we assumed that
$\tilde{f}$ maps each fiber birationally to the image.
Therefore, it suffices to show that $\pi$ is proper. 

 Since $I$ is a locally closed subscheme of
$H:=\text{Hilb}_{\tilde{M}\times\tilde{M}}(\tilde{C}\times\tilde{C})$,
which is a disjoint union of projective schemes
over $\tilde{M}\times\tilde{M}$,
and is of finite type,
if $\pi$ is not proper, then there exists a curve $S$ in $H$
which is not contained in $I$ but whose generic point is.
Now if the pullbacks of
$(\tilde{C}, \tilde{D}, \tilde{f})\times \tilde{M}$
and $\tilde{M}\times(\tilde{C}, \tilde{D}, \tilde{f})$
to the normalization $\tilde{S}$ of $S$ are isomorphic over $\tilde{S}$,
then $S$ is contained in $I$, a contradiction.
Thus it suffices to show the following: let $S$ be a smooth curve,
$P$ a point of $S$,
$p_S:C_S\rightarrow S$ and $p'_S:C'_S\rightarrow S$
families of curves induced from $C^{(d)}$
by morphisms $S\rightarrow S^{(d)}$,
$D_S$ and $D'_S$ the induced Cartier divisors
and $f_S:C_S\rightarrow X\times S$ and $f'_S:C'_S\rightarrow X\times S$
morphisms over $S$ that map any fiber of $p_S$ and $p'_S$ birationally.
Then an isomorphism of $(C_S, D_S, f_S)$ and $(C'_S, D'_S, f'_S)$
over $S\backslash \{P\}$
extends to an isomorphism over $S$. 

 If a general fiber of $C_S\rightarrow S$ is smooth,
then $C_S$ and $C'_S$ are the normalizations of
$f_S(C_S), f'_S(C'_S)\subset X\times S$,
and they coincides as they are irreducible
and have the same generic point.
By the uniqueness of normalization, $C_S$ and $C'_S$ are isomorphic
over $X\times S$.
This induces an isomorphism of $D_S$ and $D'_S$,
since it does over $S\backslash \{P\}$. 

 If a general fiber of $C_S\rightarrow S$ is singular,
then $C_S$ and $C_S'$ are isomorphic to
$(\text{reducible conic})\times S$.
As before, there exists a unique isomorphism
of the normalizations of $C_S$ and $C'_S$ over $X\times S$,
and since the restrictions of this isomorphism
to the two components of the inverse image of the double curve 
commute with the identifications by normalization maps 
(for they do over $S\backslash\{P\}$), 
it gives an isomorphism of $C_S$ and $C'_S$ over $X\times S$, 
inducing an isomorphism of $D_S$ and $D'_S$.
Thus $\pi$ is a closed immersion.

 Clearly, $I$ defines an equivalence relation.
Therefore, what we have to show is that the two projection maps
$I\rightarrow\tilde{M}$ are \'etale,
and this is equivalent to the following:
let $f_S:(C_S, D_S)\rightarrow X\times S$ be induced by
an $S$-valued point of $\tilde{M}$,
where $S$ is the spectrum
of an artinian local $\Bbb C$-algebra $(R, m)$ with $R/m=\Bbb C$,
and $P$ a geometric point of $M(H_1, H_2, H_3)$,
where $H_1, H_2, H_3$ are Cartier divisors on $X$,
corresponding to
a morphism $f_P:(C_P, D_P)\rightarrow X$
isomorphic to the restriction $f_0:(C_0, D_0)\rightarrow X$
of $f_S$ to the central fiber.
Then there exists a morphism $S\rightarrow M(H_1, H_2, H_3)$
with image $P$ which induces a family
isomorphic to $(C_S, D_S, f_S)$ over $S$, 
and for a constant family over $S=\text{Spec}\Bbb C[u]/(u^2)$
(i.e. one that is induced by a morphism to $\tilde{M}$
whose tangent is zero),
such a map is a constant map. 

 By the unique isomorphism of $C_0$ and $C_P$
commuting with $f_0$ and $f_P$,
three points $P_1, P_2, P_3$ on $C_0$ are given
satisfying $f_0^*H_i=P_i$ in a neighborhood of $P_i$.
Since $C_S$ has the same support as $C_0$,
we can take the pullback $P_{S,i}$ of $H_i$ in a neighborhood of $P_i$.
By the universal property of $M(H_1, H_2, H_3)$,
it suffices to show the following: there exists a morphism 
$g:S\rightarrow S^{(d)}$ such that there exists an isomorphism 
$h:(C_S, P_{S,1}, P_{S,2}, P_{S,3}, D_S)
\rightarrow(C^{(d)}, \sigma_1, \sigma_2, \sigma_3, D^{(d)})\times_{S^{(d)}}S$. 
Moreover, for $S=\text{Spec}\Bbb C[u]/(u^2)$ and
$(C_S, P_{S,1}, P_{S,2}, P_{S,3}, D_S)=(C_0,P_1,P_2,P_3,D_0)\times S$,
$g$ is a constant map with value $t_0$ and $h=h_0\times S$
where $h_0:(C_0, P_1, P_2, P_3, D_0)\rightarrow
(C^{(d)}_{t_0},\sigma_{1,t_0},\sigma_{2,t_0},\sigma_{3,t_0},D^{(d)}_{t_0})$ 
is the unique isomorphism.

 If $C_0$ is smooth, then $C_S\cong\Bbb P^1\times S$.
After an automorphism over $S$, we can take homogeneous coordinates $(x: y)$
such that
$P_{S,1}=(y=0)$, $P_{S,2}=(y=x)$, $P_{S,3}=(y=rx)$,
where $r\in R$, and $D_S=(x^d=0)$ hold.
The morphism $S\rightarrow S^{(d)}$ defined by $t=r$
is clearly the unique morphism
which induces $(C_S, P_{S,1}, P_{S,2}, P_{S,3}, D_S)$,
and $h$ is also unique. 

 Next, we consider the case where $C_0$ is singular.
Let $s$ be the coordinate $1/t$ of $S^{(d)}$ around $\infty$, 
as we defined before. 
Then we have the following description of $C^{(d)}$: 
$C^{(d)}=U_1\cup U_2\cup U_3$, where
$$\align{
 U_1&=\text{Spec}\Bbb C[w, s], \cr
 U_2&=\text{Spec}\Bbb C[x, z, s]/((x-s)z-s^d), \cr
 U_3&=\text{Spec}\Bbb C[y,s], \cr
}\endalign$$
and $U_1$ and $U_2$ are patched by $zw=1$ on $z\not=0, w\not=0$,
$U_2$ and $U_3$ by $xy=1$ on $x(x-s)\not=0, y(sy-1)\not=0$,
and $U_1$ and $U_3$ by $y(s^d w+s)=1$ on $s^dw+s\not=0, sy\not=0$.
$\sigma_1, \sigma_2$ are defined by $y=0, y=1$ on $U_3$,
$\sigma_3$ by $w=0$ on $U_1$
and $D^{(d)}$ by $z+\sum_{i=0}^{d-1}x^is^{d-1-i}=0$.
$(C_S, D_S)$ is obtained by replacing $\Bbb C$ by $R$
and $s$ by an element $r\in m$. 

 Consider the family
$(\tilde{C}, \tilde{D}):=(C^{(d)}, D^{(d)})\times_{S^{(d)}}\tilde{S}$ over
$\tilde{S}:=\{(s,a,b,c)\vert1/(s^dc+s),a,b\text{ are pairwise distinct 
and } c\not=-1\text{ if } d=1\}$,
and let $\tilde{P}_1, \tilde{P}_2$ be defined by $y=a, y=b$
and $\tilde{P}_3$ by $w=c$.
Then $(C_S, P_{S,1}, P_{S,2}, P_{S,3}, D_S)$ is induced
by a morphism $S\rightarrow\tilde{S}$.
Thus, to show the existence part, it suffices to show a similar statement 
for $(\tilde{C}, \tilde{P}_1, \tilde{P}_2, \tilde{P}_3, \tilde{D})$. 
We define a morphism $g:\tilde{S}\rightarrow S^{(d)}$ 
by $g^*s=(b-a)(s^dc+r)/(1-sa-s^dac)$ 
and a morphism $h:\tilde{C}\rightarrow C^{(d)}$ by 
$$\align{
 h^*s&=\frac{(b-a)(s^dc+s)}{1-sa-s^dac}, \cr
 h^*x&=\frac{(b-a)x}{1-ax}, \cr
 h^*y&=\frac{y-a}{b-a}, \cr
 h^*z&=\frac{(z-saz-s^da)(1+s^{d-1}c)^d}{1-cz}
\Bigl(\frac{b-a}{1-sa-s^dac}\Bigr)^{d-1}, \cr
 h^*w&=\frac{w-c}{(1-sa-s^daw)(1+s^{d-1}c)^d}
\Bigl(\frac{1-sa-s^dac}{b-a}\Bigr)^{d-1}. \cr
}\endalign$$
Then they give the desired cartesian diagram,
for they induce isomorphisms of fibers. 

 For the constant family over $R=\Bbb C[u]/(u^2)$,
i.e. that which is induced by $s=0$,
there is no such commutative diagram with $g^*s\not=0$:
for $d\geq2$, it is easy to see that
there are isomorphisms of $C_S$ and $C^{(d)}\times_{S^{(d)}}S$,
but they do not map $D_S$ to $D^{(d)}\times_{S^{(d)}}S$.
For $d=1$, $C_S$ and $C^{(d)}\times_{S^{(d)}}S$ are not isomorphic.
Thus $g^*s=0$, and the claim about $h$ is easy to see. 

 Thus $I$ is an \'etale equivalence relation and 
defines $\tilde{M}/I$ as a separated algebraic space. 
Patching together pieces of $\tilde{C}$, $\tilde{D}$ and $\tilde{f}$ 
by the unique isomorphism, 
we obtain a representing family. 

 The last assertion is clear from the boundedness of morphisms 
with images in a fixed homology class. 

 (2) Similar and easier. \qed
\enddemo

\proclaim{Lemma 3.2}
 Let $d$ be a positive integer 
and $P_1, P_2$ and $P_3$ different points on $\Bbb P^1$.
For Cartier divisors $H_1, H_2$ on $X$, 
let $M(H_1, H_2)$ be the scheme of morphisms $f:\Bbb P^1\rightarrow X$
with images in the class $\beta$ 
such that $f^*H_i=P_i$ in a neighborhood of $P_i$($i=1,2$) 
and that $f^*B=dP_3$ in a neighborhood of $P_3$.
Then there exists a finite collection of $H_i$'s
so that $M(H_1, H_2)$'s form an \'etale covering
of the open subspace $M_0$ of $M$ in Proposition 3.1(1)
consisting of morphisms from $\Bbb P^1$. \qed
\endproclaim

\definition{Definition 3.3}
 We denote by $M(X, B, \beta)$
the algebraic space $M$ in the Proposition 3.1 with $d=(B. \beta)$, 
and similarly for $C$, $D$, $f$. 

 Let $X$ be a smooth rational projective surface
and $B$ a smooth anti-canonical divisor of $X$.
Then a homology class of degree 2 is the same thing as a divisor class $D$,
and there are only finitely many points $P$ on $B$
with $D\vert_B\sim dP$ if $d>0$. 
So, for such a point $P$ on $B$,
we define $M(X, B, P, D)$ to be
the open and closed subspace of $M(X, B, D)$
consisting of maps whose image intersect with $B$ only at $P$ if $d>0$
and $M(X, B, P, D):=M(X, B, D)$ if $d=0$. 

 We denote their open subspaces of morphisms from $\Bbb P^1$ 
by adding $_0$. 
\enddefinition

 In the situation of the latter half of the definition, 
we may call $\deg M_0(X, B, P, D)$ 
the ``number'' of curves in $\vert D\vert$ satisfying (*) for $(B, P)$ 
by the following lemma.

\proclaim{Lemma 3.4}
 Let $X$ be a smooth rational projective surface,
$B$ a smooth anti-canonical divisor and $D$ a divisor class. Then, 

 (1) Geometric points of $M_0(X, B, P, D)$
correspond bijectively to curves $C$ on $X$ in the class $D$
which satisfy (*) for $(B, P)$. 

 (2) A point of $M_0(X, B, P, D)$ is reduced
if the normalization morphism of the corresponding curve 
is immersive outside $B$. 

 (3) $M(X, B, D)$ is supported on a finite set. 
\endproclaim
\demo{Proof}
 (1) is easy, for we required morphisms in $M(X, B, D)$ 
to be birational. 

 (3) follows from Proposition 1.1
(for morphisms from reducible conics, we look at components of them). 

 The following proposition, together with Lemma 3.2 for $d>0$, 
implies (2). 
\enddemo

\proclaim{Proposition 3.5(a log version of Proposition 3 in [Mr])}
 Let $\Cal X\rightarrow S$ be a quasi-projective family of varieties
with $S$ noetherian,
$\Cal B_i$ Cartier divisors on $\Cal X$,
$\Cal C\rightarrow S$ a flat projective family of schemes
and $\Cal D_i$ Cartier divisors on $\Cal C$ flat over $S$.
Denote the fibers by $X_s$ etc.
Let $\Cal M$ be the scheme of morphisms $f:C_s\rightarrow X_s$
such that $f^*B_{i,s}=D_{i,s}$ in a neighborhood of $\text{Supp} D_{i,s}$
and that $\sum_{i\in I}\Cal B_{i,s}$ is normal crossing over $S$
at $f(\cap_{i\in I}\text{Supp} D_{i,s})$ for any $I$:
in particular, for $I=\phi$, $\Cal X$ is smooth over $S$ at $f(C_s)$
and for $I=\{i\}$, $\Cal B_{i,s}$ is smooth over $S$
at $f(\text{Supp} D_{i,s})$.
For a point $[f:C_s\rightarrow X_s]$ of $\Cal M$,
let $\Cal T$ be the locally free sheaf on $C_s$
which is the pullback of $T_{X_s}(-\text{log}\sum_{i\in I}B_{i,s})$
near a point $P$ which is in $\cap_{i\in I}\text{Supp} D_{i,s}$
and which is not in the support of the other $D_{i,s}$'s. 
Then the completion of $\Cal M$ at $[f]$ is
defined by $h^1(C_s, \Cal T)$ equations
in the completion of $H^0(C_s, \Cal T)\times S$. 
\endproclaim
\demo{Proof}
 Similar to the proof of Proposition 3 in [Mr]. \qed
\enddemo
\noindent\qed

\definition{Definition 3.6}
 Let $l$ and $m$ be nonnegative integers,
$e$, $(a_i)_{i=1}^l$ and $(b_j)_{j=1}^m$ positive integers
with $d=d(e; (a_i); (b_j)):=3e-\sum a_i-\sum b_j\geq0$,
$B$ a smooth cubic in $\Bbb P^2$,
$P_1, \dots, P_l$ distinct points on $B$ and
$P$ a point on $B$ with $(d+\sum b_j)P+\sum a_i P_i \sim eH\vert_B$. 

 Let $(X((P_i), P, m), B((P_i), P, m; B))$,
or $(X, B)$ for short,
be the pair obtained by blowing up $P_1, \dots, P_l$
and then blowing up $m$ times
the inverse image of $P$ on the proper transform of $B$,
$E_i$ the inverse images of $P_i$
and $F_i$ the proper transform of the exceptional divisor of
the $i$-th blow-up at $P$.
In this section and the next, 
we use these notations 
and identify $B\subset\Bbb P^2$ and $B\subset X$
when it is not confusing. 

 Let $B, P_i$ and $P$ be general with
$(d+\sum b_j)P+\sum a_iP_i \sim eH\vert_B$:
we also require that
$((d+\sum b_j)/g)P+\sum (a_i/g) P_i \sim (e/g)H\vert_B$
does not hold for $g$ dividing $e$ and $a_i$. 
Let $(X, B)=(X(e;(a_i);(b_j)), B(e;(a_i);(b_j)))
:=(X((P_i), P, m), B((P_i), P, m; B))$
and $D=D(e;(a_i);(b_j)):=eH-\sum a_i E_i-\sum_{k=1}^m(\sum_{j=1}^k b_j)F_k$. 
Then we define $M(e; (a_i); (b_j))$ to be 
equal to $M_0(X, B, P, D)$ 
and we set $n(e;(a_i);(b_j)):=\deg M(e;(a_i);(b_j))$. 

 Although the variety of $(B,(P_i),P)$ satisfying the condition above 
may be reducible, 
it turns out that 
these numbers are well defined in the cases we are concerned with. 
\enddefinition

\proclaim{Lemma 3.7}
 In the situation of Definition 3.6,
if $C\in\vert D(e;(a_i);(b_j))\vert$
satisfies $C\vert_{B(e;(a_i);(b_j))}=dP$, 
it is of the form $A+\sum c_j F_j$ where $A$ is irreducible and reduced. \qed
\endproclaim

 Now we have the following relations. 

\proclaim{Theorem 3.8}
 (We use the notations in Definition 3.6.)

 (1) If $d(e; (a_i); (b_j)_{j=1}^m)>0$
and the image of any morphism representing a point 
in $M(e; (a_i); (b_j)_{j=1}^m)$ is smooth at $P$,
$n(e; (a_i); (b_j)_{j=1}^m)=n(e; (a_i); b_1,\dots,b_m,1)$. 

 (2) If $d(e; (a_i)_{i=1}^{l+1}; )=0$,
$n(e; (a_i)_{i=1}^{l+1}; )=n(e; (a_i)_{i=1}^l; a_{l+1})$. 

 (3) For $d(e; (a_i)_{i=1}^{l+1}; )>0$, assume the following: 

 (a) Any curve in $M(e; (a_i)_{i=1}^l; a_{l+1}+1)$
is transversal to $F_1$ at $P$ if $d(e;(a_i)_{i=1}^{l+1};)\mathbreak>1$, 
and is transversal to $F_1$ if $d(e;(a_i)_{i=1}^{l+1};)=1$. 

 (b) Any curve in $M(e; (a_i)_{i=1}^l; a_{l+1}+k)$
is transversal to $F_1$ for $k\geq2$. 

 Then, 
$$
 n(e; (a_i)_{i=1}^{l+1}; )=n(e; (a_i)_{i=1}^l; a_{l+1})
+(1+\delta_{d(e;(a_i)_{i=1}^{l+1};),1}a_{l+1}) n(e;(a_i)_{i=1}^l;a_{l+1}+1). 
$$
\endproclaim

\demo{Proof}
 For (1), it suffices to show the following: 
\enddemo

\proclaim{Lemma 3.9}
 Let $X$ be a smooth rational surface, 
$B\subset X$ a smooth anti-canonical divisor, 
$S$ the spectrum of an artinian local $\Bbb C$-algebra $(R, m)$
with $R/m=\Bbb C$ 
and $C_S$ a projective flat scheme over $S$
such that $C:=C_S\otimes R/m$ is a reduced curve.
Let $\sigma_S\subset C_S$ be a Cartier divisor 
mapping isomorphically onto $S$
(and hence lying on the smooth locus)
and let $D_S=d\sigma_S$.
Let $f_S:C_S\rightarrow X\times S$ be a morphism over $S$
such that $f_S^*(B\times S)=D_S$ and
that the restriction $f:=f_S\vert_C$ 
is birational onto its image and immersive at singular points. 
Then the composition of $f_S\vert_{\sigma_S}$ and the projection to $X$ 
is a constant map. 
\endproclaim

\demo{Proof}
 Let $\Cal I_S$ denote the kernel of the homomorphism
$f_S^{\#}: \Cal O_{X\times S}\rightarrow f_{S*}\Cal O_{C_S}$
and $C'_S$ the subscheme of $X\times S$ defined by $\Cal I_S$. 

 Then $\Cal I_S$ is locally free:
first, for $Q\in f(C)$ and $Q'\in f^{-1}(Q)$, 
we look for a formal function in 
$I_S:=
\ker(g_S:\hat{\Cal O}_{X\times S,Q}\rightarrow\hat{\Cal O}_{C_S,Q'})$ 
which generates 
$I:=\ker(g:\hat{\Cal O}_{X,Q}\rightarrow\hat{\Cal O}_{C,Q'})$. 

 If $Q'$ is singular, $g$ is surjective by assumption,
hence $g_S$ is also surjective. 
Then, writing down exact sequences, we see that the claim holds. 

 If $Q'$ is smooth, take coordinate functions $u\in\hat{\Cal O}_{C_S,Q'}$
and $x,y\in\hat{\Cal O}_{X\times S,Q}$ over $R$.
Replacing $u$ by another parameter in $(u)$,
we have first $y\mod m=u^a\mod m$
and then $y=u^a+\sum_{i=0}^{a-1} b_iu^i$, where $b_i\in m$.
If $x=\sum_{i=0}^{\infty} a_iu^i$,
the formal power series obtained
from $\prod_{j=1}^a \sum_{i=0}^{\infty} a_iu_j^i$
by substituting $-b_{a-1}$ for $\sum u_j$, 
$(-1)^a(b_0-y)$ for $\prod u_j$, etc., 
is an element of $I_S$ and generates $I$. 

 Taking the product over the points of $f^{-1}(Q)$,
we see that at any point $Q$
there exists a local section $\phi$ of $\Cal I_S$
which generates $\Cal I:=\ker(f^{\#}: \Cal O_X\rightarrow f_*\Cal O_C)$.
For a small neighborhood $U\subset X_S$ of $Q$,
$f_S$ is a closed immersion on $U\backslash \{Q\}$.
Since $C_S$ is flat over $S$,
$\ker(\Cal I_S\rightarrow\Cal I)=m\Cal I_S$ holds on $U\backslash \{Q\}$,
and we see that $\Cal I_S$ is generated by $\phi$ on $U\backslash \{Q\}$.
Therefore, for a section $\psi$ of $\Cal I_S$ on $U$, 
there exists a (unique) section
$\alpha\in\Gamma(U\backslash\{Q\}, \Cal O_{X\times S})$
such that $\psi=\alpha\phi$(for $\phi$ is not a zero divisor).
Then $\alpha$ extends to a section
$\tilde{\alpha}\in\Gamma(U, \Cal O_{X\times S})$ 
and $\psi=\tilde{\alpha}\phi$ holds, 
since $X$ satisfies (S$_2$). 
Thus $\Cal I_S$ is invertible. 

 Take a local analytic coordinate $u$ on $C_S$ with $\sigma_S=(u=0)$
and local parameters $(x, y)$ on $X$ such that $B=(y=0)$
and $P:=f(\sigma)=(x=y=0)$($\sigma:=\sigma_S\vert_C$). 
Then we have $f_S^{\#}y=u^d A$, 
where $A\in R[[u]]^{\times}$. 
By changing $u$ as above, we may assume $f_S^{\#}y=u^d$. 
If $f_S^{\#}x=\sum_{i=0}^{\infty}a_i u^i$, 
where $a_i\in R$ and $a_0\in m$, 
we see that $\Cal I_S\vert_B$ is generated by $(x-a_0)^d$ at $P$. 
What we have to show is that $a_0$ is zero. 

 Now $\Cal I_S^*\vert_{B\times S}$ is a line bundle on $B\times S$
with a section defining $C'_S\cap(B\times S)$,
which is given by $(x-a_0)^d=0$.
Therefore,
if we define a morphism $\alpha:S\rightarrow B$ by $x=a_0$
and $\beta:B\rightarrow\text{Pic}^d(B)$
by the line bundle $\Cal O_{B\times B}(d\Delta)$,
where $\Delta$ is the diagonal,
$\Cal I_S^*\vert_{B\times S}$ is isomorphic to
the pullback of the universal bundle
on $B\times\text{Pic}^d(B)$ by $\beta\circ\alpha$.
As $X$ is rational, there exists an isomorphism
$\Cal I_S\cong\Cal O_X(-f_*C\times S)$, 
and therefore $\beta\circ\alpha$ is a constant map.
Since $\beta$ is \'etale, $\alpha$ is also constant, 
that is, $a_0=0$. \qed
\enddemo

 (2) is clear.

 To prove (3), we look at the behavior of $M$ when $P_{l+1}$ approaches $P$. 

 Let $B'$ be a general cubic curve in $\Bbb P^2$, 
$\Cal P'_1, \dots, \Cal P'_{l+1}, \Cal P'\subset B'\times\Delta
\subset\Bbb P^2\times\Delta$
be sections over a germ $(0\in\Delta)$ of a smooth curve
and $e$ and $a_i$ positive integers
such that $\Cal P'_i$ are pairwise disjoint,
that $\Cal P'$ does not intersect $\Cal P'_i$ for $i=1, \dots, l$ 
and possibly intersects $\Cal P'_{l+1}$ only over 0
and that $d\Cal P'+\sum a_i\Cal P'_i\sim eH\vert_{B'}\times\Delta$, 
where $d=d(e; (a_i); )>0$. 

 Blow up $\Bbb P^2\times\Delta$ at $\bigcup\Cal P'_i$ 
and call this $\Cal X$,
and let $\Cal B$ be the proper transform of $B'\times\Delta$, 
$\Cal E_i$ the exceptional divisors,
$\Cal P_i$ their intersections with $\Cal B$,
$\Cal P$ the proper transform of $\Cal P'$ and
$\Cal H$ the pullback of $H\times\Delta$.
Denote their fibers over 0 by $X, B$ etc. 

 Now the construction of $M(X, B, P, D)$ can be relativized
for $(\Cal X, \Cal B, \Cal P, e\Cal H-\sum a_i\Cal E_i)\rightarrow\Delta$.
Denote the resulting algebraic space by $\Cal M$,
and the universal family by $\tilde{p}:\Cal C\rightarrow\Cal M$,
$\Cal D\subset\Cal C$ and $\tilde{f}:\Cal C\rightarrow\Cal X$.

\proclaim{Lemma 3.10}
 (1) $\Cal M$ is flat over $0\in\Delta$. 

 (2) Assume that any element $C\in\vert eH-\sum a_iE_i\vert$
with $C\vert_B=dP$ is irreducible and reduced if $P\not=P_{l+1}$, 
and is of the form $A+kE_{l+1}$ if $P=P_{l+1}$, 
where $A$ is irreducible and reduced, 
If $P=P_{l+1}$, assume furthermore that 
any curve in $M(X, B, P, eH-\sum_{i=1}^l a_iE_i-(a_{l+1}+k)E_{l+1})$
is transversal to $E_{l+1}$ for $k\geq2$. 
Then $\Cal M$ is finite over $\Delta$. 
\endproclaim

\demo{Proof}
 (1) Let $f:(C, D)\rightarrow (X, B)$
be a morphism corresponding to
a point $[f]$ of $\Cal M$.
Take general hypersurfaces $\Cal H_i$ $(i=1, 2, 3)$,
and define $\Cal N:=\Cal M(\Cal H_1, \Cal H_2, \Cal H_3)$
in the same way as $M(H_1, H_2, H_3)$.
Take a point of $\Cal N$ over $[f]$ 
representing $f:(C, D, P_1, P_2, P_3)\rightarrow(X, B, H_1, H_2, H_3)$.
Then $\Cal N$ is an \'etale neighborhood of $[f]$ in $\Cal M$. 

 Applying Proposition 3.5 to $\Cal N\rightarrow S^{(d)}\times\Delta$,
the completion of $\Cal N$ at $[f]$ is
defined by $h^1(C, \Cal T)$ equations in
the completion of $H^0(C, \Cal T)\times S^{(d)}\times\Delta$,
where $\Cal T$ is defined as in Proposition 3.5.
By Lemma 3.4, the fiber over 0 is 0-dimensional,
and by $\chi(C, \Cal T)=-1$,
we see that the equations above restrict
on $H^0(C, \Cal T)\times S^{(d)}\times\{0\}$
to the ones defining a complete intersection.
Therefore, $\Cal M$ is flat over $\Delta$.

 (2) Let $R$ be a discrete valuation ring, $K$ the quotient field of $R$,
$\text{Spec}R\rightarrow\Delta$ a morphism which maps the closed point to 0 and
$\text{Spec}K\rightarrow\Cal M$ a morphism commuting with natural morphisms.
Then what we have to show is
that there exists a finite extension $(R', K')$
such that there exists a morphism $\text{Spec}R'\rightarrow\Cal M$
commuting with natural morphisms.

 First, there exists a finite extension $K'$ of $K$
with a morphism $i:\text{Spec}K'\rightarrow\Cal N$
commuting with other morphisms,
where $\Cal N$ is an \'etale neighborhood
constructed as in the proof of (1)
for divisors $\Cal H_1, \Cal H_2, \Cal H_3$.
Since $\Cal N$ is quasi-finite over $\Delta$ 
and since the assertion is clear
if $\text{im}(i)$ is a geometric point,
we may assume that the image of $i$ is
the generic point of a curve dominating $\Delta$.
So, it suffices to show the following:
let $S$ be a smooth curve with a dominant morphism $S\rightarrow\Delta$,
$P$ a point on $S$ over 0,
$g:S_0:=S\backslash \{P\}\rightarrow\Cal N$ a morphism.
Then, by replacing $S$ by an \'etale neighborhood of $P$,
there exists a morphism $S\rightarrow\Cal M$
that commute with natural morphisms. 

 By $g$, we have a (smooth) ruled surface $C_{S_0}\rightarrow S_0$
(since the curves are irreducible by assumption),
a Cartier divisor $D_{S_0}$ whose support is a section and
a morphism $f_{S_0}: C_{S_0}\rightarrow \Cal X$
such that $f_{S_0}^*\Cal B=D_{S_0}$.
By resolution of indeterminacy and Stein factorization,
we can extend $C_{S_0}$ and $f_{S_0}$ to a morphism $p_S:C_S\rightarrow S$,
where $C_S$ is a normal surface,
and a morphism $f_S: C_S\rightarrow \Cal X$
which contracts no component of $C_P:=p_S^{-1}(P)$.
Considering $f_{S*}(C_S)\vert_{\Cal B}$,
we see that $f_{P*}(C_P)\vert_B=dP$, where $f_P:=f_S\vert_{C_P}$.
By assumption, $f_{P*}C_P=A+kE_{l+1}$,
where $A$ is irreducible, reduced and different from $E_{l+1}$.
Since $D_S:=f_S^*(\Cal B)$ is supported on a section
and does not contain any component of $C_P$,
$f_P^{-1}(B)$ consists of one point.
From this, we see that $A$ is a curve corresponding to a point
$M(X, B, P, eH-\sum_{i=1}^l a_iE_i-(a_{l+1}+k)E_{l+1})$
and that the component of $C_P$ mapping to $A$
intersects with other components at one point(if $k>0$).
By the latter fact and the assumption,
$f_P^{-1}(E_{l+1})$ contains
$a_{l+1}+k-1$ isolated points if $k>0$. 
Since $(f_S^{-1}(\Cal E_{l+1})$ is generically of degree $a_{l+1}$ 
over $\Delta$, $k$ is 0 or 1. 
Thus, $C_P$ is reduced, with one or two components, 
$f_P$ maps each fiber birationally,
and one component maps isomorphically to $E_{l+1}$ 
if $C_P$ is reducible. 

 Take general divisors $\Cal H'_1, \Cal H'_2, \Cal H'_3$.
By taking an \'etale neighborhood of $P$ in $S$,
we may assume that $p_S$ has disjoint sections $s_1, s_2, s_3$
disjoint from the singular point of $C_P$ 
that are also disjoint from $D_S$
such that the pullback of $\Cal H'_i$ is equal to $s_i$
in a neighborhood of $s_i$.
When $C_P$ is reducible,
we may assume
that $s_3$ is on the component $C_2$ mapped to $E_{l+1}$
and that $s_1$ and $s_2$ are on the other component $C_1$. 

 We show that there is a morphism $S\rightarrow S^{(d)}$
which induces $(C_S, s_1, s_2, s_3, D_S)$. 

 If $C_P$ is smooth, the assertion is clear.
Assume that $C_P$ is reducible. 
Now $C_2$ can be contracted, and if we denote the contraction
by $h:C_S\rightarrow C'_S$, $C'_S$ is smooth.
Let $s'_i=h(s_i) (i=1, 2, 3)$ and $s'_4=h(D_S)_{\text{red}}$.
Now it is easy to see the (shortest) resolution
of the indeterminacy of $h$ by blowing up smooth points
is given by blowing up $d(s'_3. s'_4)$ times the intersection of
the proper transform of $s'_3$ and the fiber over $P$.
Thus, $(C_S, s_1, s_2, s_3, D_S)$ is isomorphic to the base change of
$(C^{(d)}, \sigma_1, \sigma_2, \sigma_3, D^{(d)})$
by the morphism $S\rightarrow\Bbb P^1$ 
defined by $t=(s'_2-s'_4)(s'_3-s'_1)/(s'_2-s'_1)(s'_3-s'_4)$. 

 Finally, we obtain a morphism from $S$ 
to $\Cal M(\Cal H'_1, \Cal H'_2, \Cal H'_3)$ by the universality, 
and we are done. \qed
\enddemo

 Taking $P_i (i=1,\dots, l)$ to be general and $P_{l+1}=P$, 
the following lemma proves the theorem. 

\proclaim{Lemma 3.11}
 For $e, a_1,\dots,a_l, b$,
let the assumption be as in Theorem 3.8 (3) with $a_{l+1}=b$, 
and let $(X,B)=(X(e;(a_i);b),B(e;(a_i);b))$, $F_1$, etc.,
be as in Definition 3.6.
Let $M'$ be equal to $M_0(X, B, P, eH-\sum a_iE_i-(b+1)F_1)$
if $d:=d(e; (a_i);b)>1$
and the scheme of incidence in 
$M_0(X, B, P, eH-\sum a_i E_i-(b+1)F_1)\times F_1$
if $d=1$.
Then there is an isomorphism
$$
  j:M'\rightarrow M'':=
  M(X,B,P,eH-\sum a_iE_i-bF_1)\backslash M_0(X,B,P,eH-\sum a_iE_i-bF_1). 
$$
\endproclaim 

\demo{Proof}
 We construct $j$ as follows. 

 Case $d>1$.
Let $f_S:(C_S, D_S)\rightarrow (X, B)$ be a morphism
corresponding to a connected component $S$ of $M'$,
which is the spectrum of an artinian local $\Bbb C$-algebra.
Then we construct a morphism from $C^{(d)}_{\infty}\times S$ to $X$
by mapping the first component $C_1\times S$
by an isomorphism of 
$(C_1\times S, (D^{(d)}_{\infty}\vert_{C_1})\times S)$ and $(C_S, D_S)$ 
and the second component $C_2\times S$ by an isomorphism to $F_1\times S$
that maps $(C_1\cap C_2)\times S$ to $P$: by Lemma 3.9, they can be patched. 
For some choice of the latter isomorphism, 
the pullback of $B$ is $D^{(d)}_{\infty}\times S$, 
and this gives a morphism $S\rightarrow M''$. 

 Case $d=1$.
Let $(f_S:C_S\rightarrow (X, B), Q_S\subset F_1\times S)$
be a pair corresponding to a connected component $S$ of $M'$.
Then we construct a morphism from $C^{(1)}_{\infty}\times S$ to $X$
by mapping the first component $C_1\times S$
by an isomorphism of $C_1$ and $C_S$
so that $(C_1\cap C_2)\times S$ is mapped to $Q_S$
and the second component $C_2\times S$ by an isomorphism to $F_1\times S$
that maps $(C_1\cap C_2)\times S$ to $Q_S$
and $D^{(1)}_{\infty}$ to $B$. 
This gives a point of $M''$. 

 This map $j$ is bijective by Lemma 3.7. 

 Next, we prove the following:
let $S$ be the spectrum of an artinian local $\Bbb C$-algebra $(R, m)$
and let $(C_S, D_S)$ be the family induced by
the morphism $S\rightarrow S^{(d)}; s=\epsilon$, where $s=1/t$,
for some $\epsilon\in m$.
If there exists a morphism $f_S:C_S\rightarrow X$ such that $f_S^*B=D_S$,
then $\epsilon=0$. 

 We may assume that $\epsilon m=0$.
$C_S$ is described by replacing $\Bbb C$ by $R$ and
substituting $\epsilon$ for $s$
in the description of $C^{(d)}$ in the proof of Proposition 3.1.
Denote the central fiber by $C_0$,
the restriction of $f_S$ by $f_0$, 
$(z=0)\subset C_0$ by $A_0$, $(x=0)\subset C_0$ by $B_0$
and the point $x=z=0$ by $Q$.

 For $d\geq2$, 
this is isomorphic to $S$ times reducible conic: in fact, 
by the coordinate change $x'=x-\epsilon$ and $y'=y+\epsilon y^2$,
$C_S$ is given by $C_S=U_1\cup U_2\cup U_3$, where
$$\align{
 U_1&=\text{Spec}R[w], \cr
 U_2&=\text{Spec}R[x', z]/(x'z), \cr
 U_3&=\text{Spec}R[y'], \cr
}\endalign$$
and $U_1$ and $U_2$ are patched by $zw=1$
and $U_2$ and $U_3$ by $x'y'=1$.
Let $A_S$ denote the closed subscheme defined by $z=0$, 
$B_S$ the one defined by $x'=0$, 
each isomorphic to $\Bbb P^1_S$.
$D_S\vert_{A_S}$ is defined by $x^{\prime d-1}+d\epsilon x^{\prime d-2}=0$.
If $d>2$, we have $f_0(B_0)=F_1$.
If $d=2$, $D_S\vert_{B_S}$ is defined by $z+2\epsilon=0$ 
and we may assume $f_0(B_0)=F_1$ by interchanging $A_S$ and $B_S$. 

 Take functions $u$ and $v$ on a neighborhood of $P$ in $X$
such that $B=(v=0)$ and $F_1=(u=0)$ hold.
Since $F_1$ is a $(-1)$-curve,
$f_S\vert_{B_S}$ factors through $F_1$,
and therefore $f_S^\#u$ divides $x'$.
As $f_0(A_0)$ is transversal to $F_1$ at $P$,
we have $(f_S\vert_{A_S})^\#u\in x'\Cal O_{A_S,Q}^*$.
For $\tilde{x}:=x'+(d/(d-1))\epsilon=x-(1/(d-1))\epsilon$
and $\tilde{y}:=1/\tilde{x}$,
$(f_S\vert_{A_S})^\#v$ divides $\tilde{x}^{d-1}$
and $(f_S\vert_{A_S})^\#u\in (\tilde{x}-(d/(d-1))\epsilon)\Cal O_{A_S,Q}^*$.
By Lemma 3.9, we have $\epsilon=0$.

 For $d=1$, we can describe $C_S$ as $U_1\cup U_2\cup U_3$, where
$$\align{
 U_1&=\text{Spec}R[w], \cr
 U_2&=\text{Spec}R[x, z]/(xz-\epsilon), \cr
 U_3&=\text{Spec}R[y], \cr
}\endalign$$
and $U_1$ and $U_2$ are patched by $zw=1$,
$U_2$ and $U_3$ by $xy=1$, 
and $D_S$ as $w=0$.
We have $x\vert_{U_1\cap U_2}=\epsilon w$. 

 Let $g:X\rightarrow X'$ be the contraction of $F_1$,
and take functions $u, v$ on a neighborhood of $g(F_1)$
such that $g(B)=(v=0)$
and $(u=0)$ is a smooth curve
tangent to the branch of $g(f_0(C_0))$ corresponding to $Q$.
They define functions $u$ and $v$ in a neighborhood of $F_1$ in $X$,
$u_1$ near $F_1\backslash B$
and $v_1$ near $F_1\backslash\{f_0(Q)\}$ 
with $u=u_1 v$ and $v=v_1 u$. 

 Since $F_1$ is a $(-1)$-curve,
the restriction of $f_S$ to the closed subscheme
defined by $x=\epsilon=0$ factors through $F_1$.
Therefore, we have $(f_S\vert_{U_1})^\#u\in(\epsilon)$
and $(f_S\vert_{U'_2})^\#v\in(x,\epsilon)$,
where $U'_2=U_2\backslash\{\text{some points on }U_3\}$.
By $f_S^*B=D_S$ and Lemma 3.9, we have $(f_S\vert_{U_1})^\#u\in(w)$,
and we have $(f_S\vert_{U_1})^\#u=\epsilon wa(w)$ with $a(w)\in \Bbb C[w]$.
Since $f_0\vert_{A_0}$ is transversal to $F_1$, we have
$(f_0\vert_{A_0})^\#v\in x\Cal O_{A_0,Q}^*$,
and by $(f_S\vert_{U'_2})^\#v\in(x,\epsilon)$,
we have $(f_S\vert_{U'_2})^\#v=(xc(x,z)+\epsilon d(x,z))/(1+xb(x)) 
\mod(xz-\epsilon)$
with $c(x,z), d(x,z)\in R[x,z]$, $c(0,0)\mod m\not=0$
and $b(x)\in\Bbb C[x]$.
We may assume $c(0,0)\mod m=1$ and we have
$(f_S\vert_{U'_2})^\#v=(x+x^2c_1(x)+xc_2(x)+
\epsilon(d_0+xd_1(x)+zd_2(z)))/(1+xb(x))\mod(xz-\epsilon)$
with $d_0\in\Bbb C$, $b(x), c_1(x), d_1(x)\in\Bbb C[x]$,
$d_2(z)\in\Bbb C[z]$ and $c_2(x)\in mR[x]$.
On the other hand, we may assume that $(f_S\vert_{U_1})^\#v_1=w+e(w)$
with $e(w)\in wmR[w]$.
Now $uv_1=v$ reads 
$$
 \epsilon w^2a(w)=\epsilon(w+d_0+w^{-1}d_2(w^{-1})) 
$$
on $U_1$, 
which implies $\epsilon=0$. 

 Now we can construct $k_S:S\rightarrow M'$
for a connected component $S$ of $M''$
by restricting $f_S$ to $A_S$, the component of $C_S$
which is not mapped to $F_1$. 
Then this defines $k:M''\rightarrow M'$ 
and $j$ and $k$ are inverse to each other. \qed
\enddemo

\remark{Remark}
 From the proof(of (1) and (3)) and Proposition 2.3,
we see that, for a general point $P$ of $B$,
there exists a curve $C$ through $P$ 
which is immersed and is of degree $e$ 
such that the normalization of $C\backslash B$
is isomorphic to $\Bbb C^{\times}$. 
\endremark

\head \S4 Applications of Theorem 3.8 \endhead

 Now we give some results as corollaries of Theorem 3.8. 

\proclaim{Lemma 4.1([Kl])}
 For positive integers $e, a_i$ with $d(e; (a_i);)=1$,
$n(e; (a_i);)$ is equal to the number of rational curves 
in the blow up of $\Bbb P^2$ at general points $P_i$
in the divisor class $eH-\sum a_iE_i$,
where $E_i$ are inverse images of $P_i$. 
\endproclaim
\demo{Proof}
 If $P_i, P$ are as in Definition 3.6,
then $M_0(X((P_i,P,0), B((P_i),P,0;B), P, eH \mathbreak -\sum a_iE_i)$
is the same thing as the space of rational curves
in $X((P_i),P,0)$ in the divisor class $eH-\sum a_iE_i$.
Now move $P_i$ to general points in $\Bbb P^2$
over a parameter space $\Delta$.
Then the space of rational curves as above is flat, 
by an argument similar to that in the proof of Lemma 3.10, 
and proper, since the limit is irreducible and reduced by Lemma 3.7. \qed
\enddemo

\proclaim{Corollary 4.2}
 (1) For $e\leq4$, any curve in $M(e; (a_i); (b_j))$
has only nodes and transversal to $\sum F_j$. 

 (2) Non-zero $n(e; (a_i); (b_j))$ for $e\leq4$
with $(b_j)\not=(1,1,\dots)$ are as follows: 

The followings are 1:
$$\vbox{\halign{
  \hfil#\hfil \cr
  $n(1;1^i;)$  $(i\leq3)$, \cr
  $n(2;1^i;)$  $(i\leq5)$, \cr
  $n(3; 2, 1^i;), n(3; 1^i; 2)$  $(i\leq7)$, \cr
  $n(4;3,1^i;), n(4;1^i;3)$  $(i\leq9)$, \cr
  $n(4;2^3,1^i), n(4;2^2,1^i;2), n(4;2,1^i;2^2), n(4;1^i;2^3)$  $(i\leq6)$ \cr
}}$$

$$\vbox{\halign{
  \hfil# & \hfil# & \hfil# & \hfil# & \hfil# \cr
  $i$ & $n(3; 1^i;)$ & $n(4;2^2,1^i;)$ & $n(4;2,1^i;2)$ & $n(4; 1^i; 2^2)$\cr
  \noalign{\hrule} \noalign{\smallskip}
  0 & 3  & 4  & 4  & 2  \cr
  1 & 4  & 5  & 5  & 3  \cr
  2 & 5  & 6  & 6  & 4  \cr
  3 & 6  & 7  & 7  & 5  \cr
  4 & 7  & 8  & 8  & 6  \cr
  5 & 8  & 9  & 9  & 7  \cr
  6 & 9  & 10 & 10 & 8  \cr
  7 & 10 & 12 & 12 & 10 \cr
  8 & 12 & 12 & 12 & 10 \cr
  9 & 12 &    &    &    \cr
}}$$

$$\vbox{\halign{
  \hfil# & \hfil# & \hfil# & \hfil# \cr
  $i$ & $n(4; 2, 1^i; )$ & $n(4; 1^i; 2)$ & $n(4; 1^i; )$ \cr
  \noalign{\hrule} \noalign{\smallskip}
  0 & 11 & 10 & 16  \cr
  1 & 15 & 14 & 26  \cr
  2 & 20 & 19 & 40  \cr
  3 & 26 & 25 & 59  \cr
  4 & 33 & 32 & 84  \cr
  5 & 41 & 40 & 116 \cr
  6 & 50 & 49 & 156 \cr
  7 & 60 & 59 & 205 \cr
  8 & 72 & 71 & 264 \cr
  9 & 96 & 93 & 335 \cr
 10 & 96 & 96 & 428 \cr
 11 &    &    & 620 \cr
 12 &    &    & 620 \cr
}}$$
(Thus we have $n(4;;)=16$ again.)
\endproclaim

\demo{Proof}
 (1) We prove the assertion for $e=4$.
The cases $e\leq3$ are similar.

 (a) Any curve in $M(4;1^i;3)$ is smooth and is transversal to $F_1$. 

 Let $C$ be the curve corresponding to the unique point of $M(4;;3)$.
It is clearly smooth.
$C\cap F_1$ cannot be one point 
since $C$ and $B$ intersect at $P$ with multiplicity 9.
Let $C'$ denote the image of $C$ in $\Bbb P^2$.
The possibility left is that $C'$ has two analytic branches at $P$,
one smooth, intersecting with $B$ with multiplicity 10,
and the other an ordinary cusp.
Then we see that
$C'$ can be parameterized as $(x, y)=(u^2/(1+au^3+u^4), t=u/(1+au^3+u^4))$,
where $(x, y)$ is a coordinate on $\Bbb P^2$,
by looking at the pullbacks of two tangent lines at $P$.
Now substituting them in the equation of $B$,
we have $a=0$ and $B:y^2=-x^3+x$, which contradicts the assumption
that $P$ is of order 12. 

 For $0<i\leq9$,
consider a family whose generic fiber is $X(4;1^i;3)$
and whose central fiber is $X(4;1^{i-1};3,1)$,
where $P$ is blown up first, and then $P_i$.
We have a Cartier divisor $\Cal F_1$ whose generic fiber is $F_1$
and whose central fiber is $F_1+F_2$.
Then the limit $C_0$ of the curve $C$
representing the unique point in $M(4;1^i;3)$
lies in $M(4;1^{i-1};3,1)$ by Lemma 3.7,
and this comes from the curve $C'_0$
representing the unique point of $M(4;1^{i-1};3)$.
Since $C$ is clearly smooth,
the transversality is equivalent to reducedness of the pullback of $F_1$,
and this holds as the pullback of $F_1+F_2$ to $C_0$,
which is the same thing as the pullback of $F_1$ to $C'_0$,
is reduced by induction.

 (b) Any curve in $M(4;1^i;2^3)$ is smooth
and is transversal to $F_1+F_2+F_3$. 

 Let $C$ be the curve corresponding to the unique point of $M(4;;2^3)$.
$C$ is clearly smooth.
Since $C$ and $B$ intersect at $P$ with multiplicity 6, 
$C$ intersects with $F_3$ in two points, 
which shows the assertion for $i=0$. 

 For $0<i\leq6$, we have a family
in which the limit of the curve representing the point in $M(4;1^i;2^3)$
lies in $M(4;1^{i-1};2^3,1)$, and we have the assertion inductively
as in (a).

 (c) Any curve in $M(4;1^i;2^2)$ has one node
and is transversal to $F_1+F_2$. 

 Let $C$ be the curve corresponding to a point of $M(4;;2^2)$.
Then, since $C$ and $B$ intersect with multiplicity 8 at $P$,
$C$ is transversal to $F_2$, and it is disjoint from $F_1$. 
Assume that $C$ has a cusp.
Let $C'$ be the image of $C$ in $\Bbb P^2$
and take the quadratic transformation
defined by blowing up $P$ (on the proper transform of $B$) three times
and then blowing down the proper transform of the tangent line at $P$
and the first and second exceptional curves.
Then we have a smooth cubic $B_1$
and a cubic $C_1$ with a cusp outside of $B_1$
such that $C_1\vert_{B_1}=7P_1+2Q_1$,
where $Q_1$ is a point of order 12 on $B_1$
and $P_1$ is a point with $P_1+2Q_1\sim H\vert_{B_1}$.
Since $P_1$ is a point of order 6,
there exists a conic $D$ such that $D\vert_{B_1}=6P_1$.
Then we have $D\vert_{C_1}=6P_1$, and $P_1$ is a flex of $C_1$
for the group of smooth points on $C_1$ forms an additive group.
Thus $P_1$ is a flex of $B_1$, a contradiction.
Hence $C$ has only one node. 

 For $0<i\leq8$, we consider a smooth family $\Cal X$
over a germ $\Delta$ of smooth curve
whose general fiber is $X(4;1^i;2^2)$
and whose central fiber is $X(4;1^{i-1};2^2,1)$
with Cartier divisors $\Cal F_1$ and $\Cal F_2$
whose general fibers are $F_1$ and $F_2$
and whose central fibers are $F_1$ and $F_2+F_3$.
Changing the base if necessary,
we have a morphism from a ruled surface to $\Cal X$ over $\Delta$
whose generic fiber corresponds to a point in $M(4;1^i;2^2)$,
with image $C$.
Then the image $C_0$ of the central fiber is of the form $A+\sum a_i F_i$,
where $A$ is irreducible and reduced, by Lemma 3.7.
By $C_0\sim 4H-2F_1-4F_2-5F_3$ and $a_i\geq0$,
we see that $C_0=A$ with $[A]\in M(4;1^{i-1};2^2,1)$,
$C_0=A+F_1$ with $[A]\in M(4;1^{i-1};3,1^2)$
or $C_0=A+F_3$ with $[A]\in M(4;1^{i-1};2^3)$(only when $d\leq7$). 
We can exclude the second one noting that the inverse image of $\Cal F_1$
would contain a isolated point, by (a),
while $C\cap F_1=\phi$. 
By (b) and induction, $C$ has only nodes.
Also, the inverse image of $F_1+F_2$ to the normalization of $C$ is reduced
since the limit contains a reduced point and it is of degree 2.
For $i\leq7$, this shows that $C$ is transversal to $F_1+F_2$,
since one of the intersection is on $B$, which cannot be a node of $C$.
For $i=8$, the limit is a point of $M(4;1^7;2^2,1)$,
and we see that $C$ is transversal to $F_1+F_2$. 

 (d) Any curve in $M(4;1^i;2)$ has two nodes and is transversal to $F_1$.

 Let $C$ be the curve representing a point of $M(4;;2)$.
Then, since $C$ and $B$ intersect at $P$ with multiplicity 10,
$C$ is transversal to $F_1$.
Assume that $C$ has a cusp $Q$ and a node $R$, for example.
Let $C'$ be the image of $C$ in $\Bbb P^2$.
Denote the normalization by $f:C^\nu\rightarrow C'$ and
the curve obtained by resolving the singularities $P$ and $Q$
(resp. $P$ and $R$) of $C'$ by $C_1$(resp. $C_2$).
Take a coordinate $u$ on $C^\nu$ such that
$f^{-1}(P)=\{u=0, 1\}$, with $f^*B=11[0]+[1]$,
and that $f^{-1}(Q)=\{u=\infty\}$.
Let $f^{-1}(R)=\{u=a,b\}$.
Then, taking the line through $P$ and $R$,
we have $3([0]+[1]+[a]+[b])\sim 11[0]+[1]$ on $C_2$,
and taking the line through $P$ and $Q$
we have $3([0]+[1]+2[\infty])\sim 11[0]+[1]$ on $C_1$.
These are equivalent to $3(a+b)+2=0$ and $(a-1)^2b^8=(b-1)^2a^8$.
Thus there are finitely many possibilities for $a$ and $b$,
hence for $C$ and $P$ modulo projective equivalence.
Since $B$ and one of the branches of $C$ at $P$
intersect with multiplicity 11,
$B$ is determined uniquely by $C$ and $P$. 
Therefore, for a general $B$, there does not exist such a curve $C$.
The case when $C$ has two cusps is similar. 

 For $0<i\leq10$, we have a family such that the limit of the curve
corresponding to a point of $M(4;1^i;2)$
is $A$, $A+F_2$, $A+F_1+F_2$ or $A+F_1$,
where $A$ is a curve corresponding to a point of
$M(4;1^{i-1};2,1)$,
$M(4;1^{i-1};2,2)$($i\leq9$),
$M(4;1^{i-1};3,1)$($i\leq9$)
or $M(4;1^9;3,0)$($i=10$), respectively.
By an argument similar to that in the proof of (c),
we have the assertion inductively. 

 (e) Any curve in $M(4;1^i)$ has three nodes. 

 For $i=0$, we have shown the assertion in Theorem 2.1. 
For $0<i\leq11$, we have a family such that the limit of the curve
corresponding to a point of $M(4;1^i;)$
is $A$, $A+F_1$ or $A+2F_1$,
where $A$ is a curve corresponding to a point of
$M(4;1^{i-1};1)$,
$M(4;1^{i-1};2)$($i\leq11$)
or $M(4;1^{i-1};3)$($i\leq10$), respectively.
We can exclude the third possibility,
since the inverse image of $F_1$ would contain two isolated points by (a).
Thus $C$ has only nodes, by induction. 
The assertion for $i=12$ follows from that for $i=12$. 

 (f) The assertion for $M(4;2,1^i;)$ follows from (a) and (d),
for $M(4;2^2,1^i;)$ and $M(4;2,1^i;2)$ from (a) and (c) 
and for $M(4;2^3,1^i;)$, $M(4;2^2,1^i;2)$ and $M(4;2,1^i;2^2)$ from (b). 

 (2) $n(4;i^i;2^2)$ can be obtained as in Proposition 1.5,
and it is easy to see that 
$M(4;1^i;2^3)$ and $M(4;2,1^i;2^2)$ are 1 for $0\leq i\leq6$. 

 By Lemma 4.1, we can calculate $n(e;(a_i);)$ with $d(e;(a_i);)=1$ 
by the method of Kontsevich-Manin([CM], [KM]) and Cremona transformation. 
Then we can compute all the entry by (1) and Theorem 3.8. \qed
\enddemo

\proclaim{Corollary 4.3}
 Any quartic rational curves intersecting
with a fixed general cubic $B$ only at a fixed point $P$ of order 12
is nodal outside of $B$, and the number of such curves is as follows: 

 (a) 16: smooth at $P$, 

 (b) 10: nodal at $P$, 

 (c)  2: tacnodal at $P$, 

 (d)  1: with two smooth branches with triple contact at $P$ and 

 (e)  1: with ordinary triple point at $P$. \qed
\endproclaim
\proclaim{Corollary 4.4}
 (1) $n(5;;)=113$. \par
 (2) $n(6;;)=948$. \par
\endproclaim
\demo{Proof}
 (1) Let $P$ be a point of order 15 on $B$.
As in the proof of Theorem 2.1,
we list up the candidates for $d$ and $[a]$, with the notation there,
by requiring that $p_a\geq0$,
that the sum of two of $a_i$ does not exceed $d$
and that the sum of five of $a_i$ does not exceed $2d$.
Then we have $n(5;;)=16x+8y+z$ with 
$$\align{
  x&:=\deg M(Y, B_Y, P, 2H-E_i), \cr
  y&:=\deg M(Y, B_Y, P, 3H-\sum_{j\not=i_1,i_2}E_j), \cr
  z&:=\deg M(Y, B_Y, P, 4H-2E_i-\sum_{j\not=i}E_j), \cr
}\endalign$$
where we take $i$ such that $(2H-E_i)\vert_{B_Y}\sim5P$ for $x$, etc.,
assuming that they are well defined. 
By Lemma 2.2, Lemma 3.10 and Corollary 4.2, 
we have $z=n(4;2,1^5;)=41$, $y=n(3;1^4;)=7$ and $x=n(2;1;)=1$, 
hence the assertion. 

 (2) In the same way, we have 
$n(6;;)=21n(2;;)+27n(3;1^3;)+9n(4;2,1^4;)+3n(4;1^6;)=948$. \qed
\enddemo
\proclaim{Corollary 4.5}
 Assume that any curve in $M(e;(a_i);b)$ 
is immersed and transversal to $F_1$ for $e\leq6$. 
Then $n(7;;)=8974$. 
\endproclaim
\demo{Proof}
 We have 
$$\align{
 n(7;)=&16n(3;2;)+40n(3;1^2;)+40n(4;2,1^3;)+16n(4;1^5;)+n(5;3,1^5;) \cr
      &+8n(5;2^2,1^4;)+n(6;2^5,1;).
}\endalign$$
 By using  associativity of quantum cohomology, we have
the numbers $n(e;(a_i);)$ in Table with $d(e;(a_i);)=1$,
except for $n(6;2^8,1;)$.
Also, it is easy to see that $m(6;2^8;2)=12$.
From these, Theorem 3.8 gives Table, where $x:=n(6;2^8,1;)$.
(We obtain $n(5;;)=113$ again.) 
In particular, $n(6;;)=81n(6;2^8,1;)-6342$. 
Since this is equal to 948, we have $n(6;2^8,1;)=90$
and therefore $n(6;2^5,1;)=2419-7n(6;2^8,1;)=1789$.
Thus $n(7;;)=8974$. \qed
\enddemo

 {\bf Table.} $n(d; (a_i); b)$ for $d=5, 6$ 
under the assumption of Corollary 4.5. 

 $n(e; (a_i); b)$ is 0 or 1 for other values of $e, (a_i), b$
We set $x:=n(6;2^8,1;)$, 
and this turns out to be 90 in the proof of Corollary 4.5. 

 $d=5$. 

$$\align{
 n(5;3,2^2,1^i;)=n(5;3,2,1^i;2)=n(5;2^2,1^i;3)&=n(3;1^{i+1}), \\
 n(5;2^5,1^i;)=n(5;2^4,1^i;2)&=n(3;1^{i+4}), \\
 n(5;3,2,1^i;)=n(5;3,1^i;2)=n(5;2,1^i;3)&=n(4;2,1^i;), \\
 n(5;2^4,1^i;)&=n(4;2,1^{i+3};), \\
 n(5;2^3,1^i;2)&=n(4;1^{i+3};2), \\
 n(5;2^3,1^i;)&=n(4;1^{i+3};),
}\endalign$$

$$\vbox{\halign{
  \hfil# & \hfil# & \hfil# & \hfil# & \hfil# & \hfil# \cr
  $i$ & $n(5;2^2,1^i;2)$ & $n(5;3,1^i;)$ & $n(5;1^i;3)$
    & $n(5;2^2,1^i;)$ & $n(5;2,1^i;2)$ \cr
  \noalign{\hrule} \noalign{\smallskip}
  0 & 55  & 24  & 23  & 100  & 89   \cr
  1 & 79  & 35  & 34  & 155  & 140  \cr
  2 & 110 & 50  & 49  & 234  & 214  \cr
  3 & 149 & 70  & 69  & 344  & 318  \cr
  4 & 197 & 96  & 95  & 493  & 460  \cr
  5 & 255 & 129 & 128 & 690  & 649  \cr
  6 & 325 & 170 & 169 & 945  & 895  \cr
  7 & 416 & 220 & 219 & 1270 & 1210 \cr
  8 & 584 & 280 & 279 & 1686 & 1614 \cr
  9 & 620 & 352 & 351 & 2270 & 2174 \cr
 10 &     & 448 & 447 & 3510 & 3222 \cr
 11 &     & 640 & 636 & 3510 & 3510 \cr
 12 &     & 640 & 640 &      &      \cr
}}$$

$$\vbox{\halign{
  \hfil# & \hfil# & \hfil# & \hfil# \cr
  $i$ & $n(5;2,1^i;)$ & $n(5;1^i;2)$ & $n(5;1^i;)$ \cr
  \noalign{\hrule} \noalign{\smallskip}
  0 & 127 & 104 & 113\cr
  1 & 216 & 182 & 217\cr
  2 & 356 & 307 & 399\cr
  3 & 570 & 501 & 706\cr
  4 & 888 & 793 & 1207\cr
  5 & 1348 & 1220 & 2000\cr
  6 & 1997 & 1828 & 3220\cr
  7 & 2892 & 2673 & 5048\cr
  8 & 4102 & 3823 & 7721\cr
  9 & 5716 & 5365 & 11544\cr
  10 & 7890 & 7443 & 16909\cr
  11 & 11112 & 10476 & 24352\cr
  12 & 18132 & 16212 & 34828\cr
  13 & 18132 & 18132 & 51040\cr
  14 &       &       & 87304\cr
  15 &       &       & 87304\cr
}}$$

 $d=6$. 

$$\align{
 n(6;4,2^3,1^i;)=n(6;4,2^2,1^i;2)=n(6;2^3,1^i;4)&=n(3;1^{i+1}), \\
 n(6;3^3,1^i;)=n(6;3^2,1^i;3)&=n(3;1^i), \\
 n(6;3^2,2^3,1^i;)=n(6;3^2,2^2,1^i;2)=n(6;3,2^3,1^i;3)&=n(3;1^{i+3};), \\
 n(6;3,2^6, 1^i;)=n(6;3,2^5,1^i;2)=n(6;2^6,1^i;3)&=n(3;1^{i+6}), \\
 n(6;2^9;)=n(6;2^8;2)&=12, \\
 n(6;4,2^2,1^i;)=n(6;4,2,1^i;2)=n(6;2^2,1^i;4)&=n(4;2,1^i;), \\
 n(6;3^2,2^2,1^i;)=n(6;3,2^2,1^i;3)&=n(4;2,1^{i+2};), \\
 n(6;3^2,2,1^i;2)&=n(4;1^{i+2};2), \\
}\endalign$$
$$\align{
 n(6;3,2^5,1^i;)=n(6;2^5,1^i;3)&=n(4;2,1^{i+5};), \\
 n(6;3,2^4,1^i;2)&=n(4;1^{i+5};2), \\
 n(6;4,2,1^i;)=n(6;4,1^i;2)=n(6;2,1^i;4)&=n(5;3,1^i;), \\
 n(6;3^2,2,1^i;)=n(6;3,2,1^i;3)&=n(4;1^{i+2};), \\
 n(6;3^2,1^i;2)&=n(5;2^2,1^{i-1};2)(i>0), \\
 n(6;3^2;2)&=37, \\
 n(6;3,2^4,1^i)&=n(4;1^{i+5}), \\
}\endalign$$

$$\vbox{\halign{
  \hfil# & \hfil# & \hfil# & \hfil# & \hfil# & \hfil# \cr
  $i$ & $n(6;2^8,1^i)$ & $n(6;2^7,1^i;2)$
    & $n(6;2^7,1^i;)$ & $n(6;2^6,1^i;2)$ \cr
  \noalign{\hrule} \noalign{\smallskip}
  0 & $x-24$ & $x-25$ & $604-4x$ & $595-4x$ \cr
  1 & $x$    & $x-3$  & $579-3x$ & $569-3x$ \cr
  2 & $x$    & $x$    & $576-2x$ & $564-2x$ \cr
  3 &        &        & 576      & 540      \cr
  4 &        &        & 576      & 576      \cr
}}$$

$$\align{
 n(6;3^2,1^i;)=n(6;3,1^i;3)&=n(5;2^2,1^{i-1})(i>0), \\
 n(6;3^2;)=n(6;3;3)&=63, \\
 n(6;3,2^3,1^i;)&=n(5;2^2,1^{i+2};), \\
 n(6;3,2^2,1^i;2)&=n(5;2,1^{i+2};2), \\
}\endalign$$

$$\vbox{\halign{
  \hfil# & \hfil# & \hfil# & \hfil# & \hfil# & \hfil# \cr
  $i$ & $n(6;4,1^i)$ & $n(6;1^i;4)$
    & $n(6;2^3,1^i;3)$ & $n(6;2^6,1^i;)$ & $n(6;2^5,1^i;2)$ \cr
  \noalign{\hrule} \noalign{\smallskip}
  0 & 46   & 45   & 230  & $-180+9x$ & $-221+9x$ \cr
  1 & 70   & 69   & 339  & $415+5x$  & $365+5x$  \cr
  2 & 105  & 104  & 487  & $984+2x$  & $924+2x$  \cr
  3 & 155  & 154  & 683  & 1548      & 1476      \cr
  4 & 225  & 224  & 937  & 2088      & 1992      \cr
  5 & 321  & 320  & 1261 & 3240      & 2952      \cr
  6 & 450  & 449  & 1676 & 3240      & 3240      \cr
  7 & 620  & 619  & 2258 & & \cr
  8 & 840  & 839  & 3462 & & \cr
  9 & 1120 & 1119 & 3510 & & \cr
 10 & 1472 & 1471 &      & & \cr
 11 & 1920 & 1919 &      & & \cr
 12 & 2560 & 2559 &      & & \cr
 13 & 3840 & 3835 &      & & \cr
 14 & 3840 & 3840 &      & & \cr
}}$$

$$
 n(6;3,2^2,1^i;)=n(5;2,1^{i+2}), 
$$

$$\vbox{\halign{
  \hfil# & \hfil# & \hfil# & \hfil# & \hfil# \cr
  $i$ & $n(6;3,2,1^i;2)$ & $n(6;2^2,1^i;3)$
    & $n(6;2^5,1^i;)$ & $n(6;2^4,1^i;2)$ \cr
  \noalign{\hrule} \noalign{\smallskip}
  0 & 316   & 345   & $2640-16x$ & $2525-16x$ \cr
  1 & 511   & 555   & $2419-7x$  & $2264-7x$  \cr
  2 & 804   & 868   & $2784-2x$  & $2580-2x$  \cr
  3 & 1232  & 1322  & 3708       & 3445       \cr
  4 & 1841  & 1964  & 5184       & 4850       \cr
  5 & 2687  & 2851  & 7176       & 6749       \cr
  6 & 3838  & 4052  & 10128      & 9512       \cr
  7 & 5381  & 5656  & 16608      & 14748      \cr
  8 & 7462  & 7818  & 16608      & 16608      \cr
  9 & 10492 & 11016 &            &            \cr
 10 & 16272 & 17748 &            &            \cr
 11 & 18132 & 18132 &            &            \cr
}}$$

$$\vbox{\halign{
  \hfil# & \hfil# & \hfil# & \hfil#
    & \hfil# & \hfil# \cr
  $i$ & $n(6;3,2,1^i;)$ & $n(6;3,1^i;2)$ & $n(6;2,1^i;3)$
    & $n(6;2^4,1^i;)$ & $n(6;2^3,1^i;2)$ \cr
  \noalign{\hrule} \noalign{\smallskip}
  0 & 444   & 381   & 420   & $-473+25x$ & $-703+25x$ \cr
  1 & 760   & 660   & 725   & $2052+9x$  & $1713+9x$ \cr
  2 & 1271  & 1116  & 1221  & $4316+2x$  & $3829+2x$ \cr
  3 & 2075  & 1841  & 2005  & 6896       & 6213      \cr
  4 & 3307  & 2963  & 3211  & 10341      & 9404      \cr
  5 & 5148  & 4655  & 5019  & 15191      & 13930     \cr
  6 & 7835  & 7145  & 7665  & 21940      & 20264     \cr
  7 & 11673 & 10728 & 11453 & 31452      & 29194     \cr
  8 & 17054 & 15784 & 16774 & 46200      & 42738     \cr
  9 & 24516 & 22830 & 24164 & 79416      & 68886     \cr
 10 & 35008 & 32738 & 34560 & 79416      & 79416     \cr
 11 & 51280 & 47770 & 50640 &            &           \cr
 12 & 87544 & 77014 & 84984 &            &           \cr
 13 & 87544 & 87544 & 87544 &            &           \cr
}}$$

$$\vbox{\halign{
  \hfil# & \hfil# & \hfil#
    & \hfil# & \hfil# \cr
  $i$ & $n(6;3,1^i;)$ & $n(6;1^i;3)$
    & $n(6;2^3,1^i;)$ & $n(6;2^2,1^i;2)$ \cr
  \noalign{\hrule} \noalign{\smallskip}
  0 & 459    & 414    & $5340-36x$ & $4995-36x$ \cr
  1 & 840    & 771    & $4637-11x$ & $4082-11x$ \cr
  2 & 1500   & 1396   & $6350-2x$  & $5482-2x$  \cr
  3 & 2616   & 2462   & 10179      & 8857       \cr
  4 & 4457   & 4233   & 16392      & 14428      \cr
  5 & 7420   & 7100   & 25796      & 22945      \cr
  6 & 12075  & 11626  & 39726      & 35674      \cr
  7 & 19220  & 18601  & 59990      & 54334      \cr
  8 & 29948  & 29109  & 89184      & 81366      \cr
  9 & 45732  & 44613  & 131922     & 120906     \cr
 10 & 68562  & 67091  & 200808     & 183060     \cr
 11 & 101300 & 99381  & 359640     & 305244     \cr
 12 & 149070 & 146511 & 359640     & 359640     \cr
 13 & 226084 & 222249 &            &            \cr
 14 & 401172 & 385812 &            &            \cr
 15 & 401172 & 441072 &            &            \cr
}}$$

$$\vbox{\halign{
  \hfil# & \hfil# & \hfil#
    & \hfil# & \hfil# & \hfil# \cr
  $i$ & $n(6;2^2,1^i;)$ & $n(6;2,1^i;2)$
    & $n(6;2,1^i;)$ & $n(6;1^i;2)$ & $n(6;1^i;)$ \cr
  \noalign{\hrule} \noalign{\smallskip}
  0 & $-2381+49x$ & $-2801+49x$ & $7336-64x$ & $6922-64x$ & $-6342+81x$ \cr
  1 & $2614+13x$  & $1889+13x$  & $4535-15x$ & $3764-15x$ & $580+17x$   \cr
  2 & $6696+2x$   & $5475+2x$   & $6424-2x$  & $5028-2x$  & $4344+2x$   \cr
  3 & 12178       & 10173       & 11899      & 9437       & 9372        \cr
  4 & 21035       & 17824       & 22072      & 17839      & 18809       \cr
  5 & 35463       & 30444       & 39896      & 32796      & 36648       \cr
  6 & 58408       & 50743       & 70340      & 58714      & 69444       \cr
  7 & 94082       & 82629       & 121083     & 102482     & 128158      \cr
  8 & 148416      & 131642      & 203712     & 174603     & 230640      \cr
  9 & 229782      & 205618      & 335354     & 290741     & 405243      \cr
 10 & 350688      & 316128      & 540972     & 473881     & 695984      \cr
 11 & 533748      & 483108      & 857100     & 757719     & 1169865     \cr
 12 & 838992      & 754008      & 1340208    & 1193697    & 1927584     \cr
 13 & 1558272     & 1295640     & 2094216    & 1871967    & 3121281     \cr
 14 & 1558272     & 1558272     & 3389856    & 3004044    & 4993248     \cr
 15 &             &             & 6506400    & 5302884    & 7997292     \cr
 16 &             &             & 6506400    & 6506400    & 13300176    \cr
 17 &             &             &            &            & 26312976    \cr
 18 &             &             &            &            & 26312976    \cr
}}$$


\Refs

\widestnumber\key{KM}
\ref \key B \by F.~Beukers 
  \paper Ternary form equations
  \paperinfo Utrecht Univ. Dep. of Math. Preprint, No 771 \yr 1993
  \endref
\ref \key CM \by B.~Crauder and R.~Miranda
  \paper Quantum cohomology of rational surfaces
  \inbook The Moduli Space of Curves 
  \eds R.~Dijkgraaf, C.~Faber and G.~van der Geer
  \bookinfo Prog. in Math. \vol 129
  \publ Birkh\"auser \yr 1995 \pages 33--80
  \endref
\ref \key Kl \by J.~Koll\'ar 
  \paper Rational curves in the plane 
  \paperinfo Lecture at Warwick symposium on 3-folds, Dec. 1995
  \endref
\ref \key Kn \by M.~Kontsevich
  \paper Enumeration of rational curves via torus actions
  \inbook The Moduli Space of Curves 
  \eds R.~Dijkgraaf, C.~Faber and G.~van der Geer
  \bookinfo Prog. in Math. \vol 129
  \publ Birkh\"auser \yr 1995 \pages 335--368
  \endref 
\ref \key KM \by M.~Kontsevich and Yu.~Manin
  \paper Gromov-Witten classes, quantum cohomology and enumerative geometry
  \paperinfo MPI preprint and hep-th/9402147 \yr 1994
  \endref
\ref \key Mr \by S.~Mori 
  \paper Projective manifolds with ample tangent bundles
  \jour Ann. of Math. (2) \vol 110 \yr 1979 \pages 593--606
  \endref
\ref \key Ms \by D.~Morrison 
  \paper Mirror symmetry and rational curves on quintic threefolds: 
   a guide for mathematicians
  \jour J. Amer. Math. Soc. \vol 6 \yr 1993 \pages 223--247
  \endref 
\ref \key R \by Z.~Ran 
  \paper The number of unisecant rational cubics to a plane cubic
  \paperinfo preprint
  \yr 1995
  \endref
\ref \key RT \by Y.~Ruan and G.~Tian 
  \paper A mathematical theory of quantum cohomology \paperinfo preprint
  \yr 1993
  \endref
\ref \key V \by R.~Vidunas 
  \paper Rational curves intersecting an elliptic curve in $\Bbb P^2$ 
  at one point and application to arithmetic geometry 
  \paperinfo preprint \yr 1996 
  \endref
\ref \key X \by G.~Xu 
  \paper On the intersections of rational curves with cubic plane curves
  \paperinfo Duke e-print, alg-geom/9511006 \yr 1995 
  \endref

\endRefs
\enddocument